\documentclass{amsart}
\usepackage{amssymb}
\usepackage{amscd}
\usepackage{graphics}

\usepackage[
hyperref
]{degt}

\def\Dg:{{\bf Dg:}\enspace\ignorespaces}

\def\bL{\bold L}

\def\bS{\bold S}
\def\bT{\bold T}
\def\bX{\bold X}
\let\sset\bS
\def\smin{\sset\subtext{min}}

\def\htype{\Cal H}

\def\tS{\tilde\sset_h}

\let\Gm\mu

\let\CX\CalX
\let\CY\CalY
\let\CZ\CalZ

\def\FF#1{\Bbb F_{#1}}
\def\CG#1{\Bbb G_{#1}}

\def\Dynkin{\frak G}
\def\SSG{\Bbb S}

\def\SSGp{\SSG_+}

\def\refl#1{t_{#1}}
\def\pss{\frak o}  
\def\dso{\frak s}  

\def\dual{^\sharp}
\def\units{^\times}
\def\map{\QOPNAME{d}}
\def\map{\mathrm{d}}
\let\dN=\CalN
\def\dL{\Cal L}
\def\dS{\Cal S}
\def\dT{\Cal T}

\let\GG=\Gamma
\def\GQ{\GG_\Q}
\def\GA{\GG_{\Bbb A}}
\def\GAz{\GG_{\Bbb A,0}}

\let\one=1  
\def\quo#1{\<#1\>}
\def\HH{H_\infty}
\let\HH=H

\def\MM{\mathrm{E}}
\let\mm\Gg
\def\mma{\bar\mm}
\def\mmb{\bar\Gb}
\def\mmh{\frak m}
\def\mmh{\mathrm{e}}
\def\plus{^+}
\let\genus=g
\def\tSigma{\tilde\Sigma}
\def\theN{(\hyperlink{convention}{*})}
\def\Oplus{\OG_+}

\def\torus#1{$#1$-torus}
\def\Sym{\QOPNAME{Sym}}
\def\spin{\QOPNAME{spin}}
\def\ds{\QOPNAME{ds}}
\def\conv{\QOPNAME{conv}}
\let\fe\into            
\let\pe\rightarrowtail  
\let\dg\pe              
\let\gd\leftarrowtail
\def\CK{\Cal K}         

\def\all{\Cal M}
\let\cluster\CalC
\let\ns\varnothing
\let\ns1

\let\specialmark\star

\let\fiber=F

\let\sextic=D

\let\set\frak

\def\sC{\set C}
\def\tC{\tilde\sC}

\def\tabsize{\small}
\let\tabquad\relax
\def\tabrefform#1{$^{[#1]}$}
\def\maketabref#1{\tabrefform{\ref{#1}}}
\def\tabref{\ifx\REF\empty\else\smash{\llap{\maketabref{\REF}\,}}\fi}
\def\tabbox#1{\hbox to.31\hsize{#1\hss}}
\def\tabentry#1//#2\end{\tabbox{\tabquad\xdef\REF{#2}\tabref\singset{#1}}}
\def\tabboxii#1{\hbox to.28\hsize{#1\hss}}
\def\tabentryii#1\end{\tabboxii{\singset{#1}}}
\def\tab{\vtop\bgroup\def\tabbreak{\endtab\qquad\tab}%
 \def\endtab{\noalign{\vspace{2pt}\hrule}\egroup\egroup}%
 \tabsize\let\\\cr\catcode`\|\the\catcode`\&%
 \halign\bgroup\quad\strut\tabentry##\end\hss\quad&$##$\quad\cr
 \noalign{\hrule\vspace{2pt}}%
 \omit\strut\quad Singularities\hss&(r,c)\cr
 \noalign{\vspace{1pt}\hrule\vspace{2pt}}}
\def\tabii{\vtop\bgroup\def\tabbreak{\endtab\qquad\tabii}%
 \def\endtab{\egroup\egroup}%
 \tabsize\halign\bgroup\strut\tabentryii##\end\hss\cr}
\def\tabs{\hrule height0pt\hbox to\hsize\bgroup\hss}
\def\endtabs{\hss\egroup}

\theoremstyle{definition}
 \addtheorem{Notation}{notation}

\title{Geography of irreducible plane sextics}

\author{Ay\c{s}eg\"{u}l Akyol}

\address{%
Abdullah G\"{u}l University\\
Faculty of Engineering and Natural Sciences\\
38039 Kayseri, Turkey}

\email{aysegul.akyol@agu.edu.tr}

\author{Alex Degtyarev}

\address{%
Department of Mathematics\\
Bilkent University\\
06800 Ankara, Turkey}

\email{degt@fen.bilkent.edu.tr}



\begin{document}

\begin{abstract}
We complete the equisingular deformation classification of
irreducible singular plane
sextic curves. As a by-product, we
also compute the fundamental groups of the
complement of all but a few maximizing sextics.
\end{abstract}

\maketitle

\section{Introduction}

Throughout the paper, all
varieties are over the field~$\C$ of complex numbers.

Our principal result is the completion of the classification of irreducible
plane sextics (curves of degree~$6$) up to equisingular deformation.
We confine ourselves to \emph{simple} sextics only, \ie, those with
$\bA$--$\bD$--$\bE$ singularities (see \autoref{s.sextics}). The non-simple
ones require completely different techniques
and are well known; surprisingly, their study is much easier:
the statements were announced by the second author long
ago, and formal proofs can be found in~\cite{degt:book}.
Note also that degree~$6$ is the first nontrivial case (see~\cite{degt:book}
for the statements on quintics, and
quartics were already known to Klein;
see also M.~Namba~\cite{Namba} for an excellent account of the sets of singularities
realized by curves of degree up to five)
and, probably, the last case that can be
completely understood, thanks to the close relation between plane sextics and
$K3$-surfaces.

The systematic study of simple sextics based on the theory of $K3$-surfaces
was initiated by U.~Persson~\cite{Persson:sextics}, who proved
that the total Milnor number~$\Gm$ of such a curve does not
exceed~$19$.
Based on this approach,
T.~Urabe~\cite{Urabe:sextics} listed the possible sets of singularities with
$\Gm\le16$, and this result was extended to a complete list
of the sets of singularities realized by simple sextics
by J.~G.~Yang~\cite{Yang}. Later, using the arithmetical
reduction~\cite{degt:JAG},
I.~Shimada~\cite{Shimada:maximal} gave a complete description of the
moduli
spaces of the \emph{maximizing} ($\Gm=19$) sextics.
In the meanwhile, a number of independent
(not explicitly related to the $K3$-surfaces)
attempts to attack the classification problem has also been made, see, \eg,
\cite{Artal:trends,Artal:invariants} (defining equations of a number of
maximizing sextics),
\cite{Oka.Pho:moduli,Oka.Pho:groups} (sets of singularities and
explicit equations of sextics of torus type),
\cite{degt:symmetric,degt:8a2,degt:Oka3} (sextics admitting stable
projective symmetries),
\cite{degt:book} (sextics with a triple point),
\etc.


At some point it was clearly understood,
partially in conjunction with Oka's
conjecture~\cite{Oka:conjecture} and partially due to the arithmetical
reduction of the problem~\cite{degt:JAG},
that irreducible sextics~$\sextic$ should be subdivided into classes
according to the maximal generalized dihedral
quotient~$Q_\sextic$ that the fundamental
group $\pi_1(\Cp2\sminus\sextic)$ admits.
If this quotient is large, $\ls|Q_\sextic|>6$, the curves are relatively few
in number and can easily be listed manually (see~\cite{degt:Oka} and
\autoref{s.special}),
using Nikulin's sufficient
uniqueness conditions~\cite{Nikulin:forms}.
The present paper fills the gap and covers the two remaining cases:
non-special sextics ($Q_\sextic=0$, see \autoref{th.classification}) and
$1$-torus sextics ($Q_\sextic=\DG6$, see \autoref{th.torus}).
On the arithmetical side,
our computation is based on the stronger (non-)uniqueness
criteria due to
Miranda--Morrison
\cite{Miranda.Morrison:1,Miranda.Morrison:2,Miranda.Morrison:book}.
For
an even
further illustration of the power of
\cite{
Miranda.Morrison:book},
we solve
a few more subtle geometric problems, namely,
we compute
the monodromy representation of the fundamental groups of the
equisingular strata (in
other words, we classify sextics with marked
singular points, see \autoref{s.permutations} and \autoref{th.group}),
we discuss whether the strata are real and whether they contain real curves
(the interesting discovery here is \autoref{prop.A7+A6+A5}),
and we give a complete description of the adjacencies of the strata
(see \autoref{s.adjacencies} and Propositions~\ref{prop.p=2},
\ref{prop.p=3}, \ref{prop.p=7}).

There are three sets of singularities
that deserve special attention: to the best of our
knowledge, phenomena of this kind have not been observed before.
It is quite common that the (discrete) moduli spaces of maximizing sextics
are disconnected, see~\cite{Shimada:maximal}. For about a dozen of the sets
of singularities with $\Gm=18$, the moduli space (of dimension~$1$) consists
of two complex conjugate components (see \autoref{tab.disconnected}; the
first such example,
\viz. \singset{E6+A11+A1}, was found in~\cite{Aysegul:paper}).
We discover a set of
singularities, \viz. \singset{E6+2A5+A1}, with $\Gm=17$ and disconnected
moduli space (two conjugate components of dimension~$2$),
and another one, \singset{2A9}, with
$\Gm=18$ and the moduli space consisting of two disjoint \emph{real}
components (see \autoref{prop.2A9}).
Finally, the moduli space
corresponding to
the set of singularities
\singset{A7+A6+A5}, $\Gm=18$, consists of a
single component, which is hence real, but it contains no real curves (see
\autoref{prop.A7+A6+A5}).

As another important by-product of
Theorems~\ref{th.classification}
and~\ref{th.torus},
we obtain Corollaries \ref{cor.pi1}
and~\ref{cor.pi1.torus}, computing the fundamental groups of the complements
of all but a few maximizing irreducible sextics.
In fact, no computation is found in this paper: we merely use the
classification, the degeneration principle, and previously known groups.
Most statements on the fundamental
groups were known conjecturally;
more precisely, the groups of \emph{some} sextics with
certain sets of singularities were known, and our principal
contribution is the connectedness of the moduli spaces.

\subsection{Contents of the paper}
The principal results of the paper are stated in \autoref{S.results},
after the necessary terminology and notation have been introduced.
For the reader's convenience, we
also discuss the other irreducible simple sextics (see \autoref{s.special})
and list the known fundamental groups.
In \autoref{S.lattices}, we recall the fundamentals of Nikulin's theory of
discriminant forms and lattice extensions, give a brief introduction to
Miranda--Morrison's theory~\cite{Miranda.Morrison:book}, and recast some of
their results in a form more suitable for our computations.
In \autoref{S.classification}, we recall the notion of
(abstract) homological type and
the arithmetical reduction~\cite{degt:JAG} of the classification problem
(see \autoref{s.htype} and \autoref{s.steps})
and begin the proof of our
principal results, classifying the
plane sextics up to equisingular deformation
\emph{and}
complex conjugation.
As a digression, we classify also sextics with marked singular points, see
\autoref{s.permutations}.
With the classification in hand,
the computation of the fundamental groups is almost straightforward;
it is outlined in \autoref{S.pi1}.
Finally, in \autoref{S.real}, we discuss real strata and real curves,
completing the deformation classification of simple sextics.
As another digression, in \autoref{s.adjacencies} we describe the adjacencies
of the non-real strata.
A few further results obtained after this paper was submitted are outlined
briefly
in \autoref{s.reducible}.

\subsection{Acknowledgements}
We are grateful to V.~Nikulin, who drew our attention to Miranda--Morrison's
works~\cite{Miranda.Morrison:1,Miranda.Morrison:2}.
To large extent, this text was written during the second author's stay at
the
\emph{Max-Planck-Institut f\"{u}r Mathematik}, Bonn, partially supported by
the ``Tropical Geometry and Topology'' program.
We
would like to
extend our gratitude
to the institute and its friendly staff and to the organizers of
the program.

\section{Principal results}\label{S.results}

\subsection{Notation}\label{s.notation}
We use the notation $\CG{n}:=\Z/n\Z$ (reserving $\Z_p$ and $\Q_p$ for
$p$-adic numbers) and $\DG{2n}$ for the cyclic group
of order~$n$ and dihedral group
of order $2n$, respectively.
As usual, $\SL(n,\Bbbk)$ is the group of $(n\times n)$-matrices~$M$ over a
field~$\Bbbk$ such that $\det M=1$.

The notation $\BG{n}$ stands for the braid group on $n$ stings.
The \emph{reduced braid group} (or the \emph{modular group}) is the quotient
$\MG=\BG3/(\Gs_1\Gs_2)^3$ of $\BG3$ by its center;
one has $\MG=\PSL(2,\Z)=\CG2*\CG3$.
The braid group is generated by the \emph{Artin generators}
$\sigma_i$, $i=1,\ldots,n-1$,
subject to the relations
\[*
[\Gs_i, \Gs_j]=1\quad\text{if $\ls|i-j|>1$}, \qquad
\Gs_i \Gs_{i+1} \Gs_i = \Gs_{i+1} \Gs_i \Gs_{i+1}.
\]

Throughout the paper, all group actions are right, and we use the notation
$(x,g)\mapsto x\ra g$.
The standard action of $\BG{n}$ on the free group $\<\Ga_1,\ldots,\Ga_n\>$
is as follows:
\[*
\Gs_i:\begin{cases}
    \Ga_i\mapsto \Ga_i\Ga_{i+1}\Ga_i^{-1}, & \hbox{ } \\
    \Ga_{i+1}\mapsto\Ga_i, & \hbox{ } \\
    \Ga_j\mapsto\Ga_j, & \hbox{if $j\neq i,i+1$}
  \end{cases}
\]
The element $\Gr_n:=\Ga_1\ldots\Ga_n\in\<\Ga_1,\ldots,\Ga_n\>$ is
preserved by~$\BG{n}$. Given a pair $\Ga_1, \Ga_2$, we use the notation
$\{\Ga_1,\Ga_2\}_n:=\Ga_2\1(\Ga_2\ra\Gs_1^n)\in\<\Ga_1,\Ga_2\>$ for $n\in\Z$.
Explicitly, the relation $\{\Ga_1,\Ga_2\}_n=1$ in a group boils down to
\[*
\aligned
(\Ga_1\Ga_2)^k=(\Ga_2\Ga_1)^k,&\quad\text{if $n=2k$ is even},\\
(\Ga_1\Ga_2)^k\Ga_1=(\Ga_2\Ga_1)^k\Ga_2,&\quad\text{if $n=2k+1$ is odd}.
\endaligned
\]
In particular, $\{\Ga_1,\Ga_2\}_1=1$ means $\Ga_1=\Ga_2$, and
$\{\Ga_1,\Ga_2\}_2=1$ means $[\Ga_1,\Ga_2]=1$.

We denote by $\PP=\{2,3,\ldots\}$ the set of all primes.

The
group of units of a commutative ring~$R$ is denoted by $R\units$.
We recall
that
$\Z_p\units/(\Z_p\units)^2=\{\pm1\}$ for $p\in\PP$ odd,
and $\Z_2\units/(\Z_2\units)^2=(\Z/8)\units\cong\{\pm1\}\times\{\pm1\}$ is
generated by~$7\bmod8$ and~$5\bmod8$.
If $m\in\Z$ is prime
to~$p$, its class in $\Z_p\units/(\Z_p\units)^2$ is
the Legendre symbol $(\frac{m}p)\in\{\pm1\}$ if $p$ is odd or
$m\bmod8\in(\Z/8)\units$ if $p=2$.

\subsection{Simple sextics}\label{s.sextics}
A \emph{sextic} is a plane curve $\sextic\subset\Cp2$ of degree six.
A sextic is \emph{simple} if all its singular points are simple, \ie, those
of type $\bA$--$\bD$--$\bE$, see~\cite{Durfee}.
If this is the case, the minimal resolution of
singularities~$X$ of the double covering of $\PP^2$ ramified at $\sextic$ is
a $K3$-surface.
The intersection index form
$H_2(X)\cong2\bE_8 \oplus 3\bU$ is (the only)
even unimodular lattice of signature $(\sigma_+,\sigma_-)=(3,19)$
(see \autoref{s.extensions};
here, $\bU$ is the hyperbolic plane).
We fix the notation $\bL:=2\bE_8 \oplus 3\bU$.

For each simple singular point~$P$ of $\sextic$,
the components
of the exceptional divisor $E\subset X$ over~$P$ span
a root lattice in~$\bL$ (see \autoref{s.root}).
The (obviously orthogonal) sum of these sublattices is denoted by
$\sset(\sextic)$ and is referred to as the \emph{set of singularities}
of~$\sextic$.
(Recall that the types of the individual singular points are uniquely
recovered from $\sset(\sextic)$, see \autoref{s.root}.)
The rank $\rank\sset(\sextic)$ equals the total Milnor number $\Gm(\sextic)$.
Since $\sset(\sextic)\subset\bL$ is
negative definite,
one has $\Gm(\sextic)\le19$, see \cite{Persson:sextics}.
If $\Gm(\sextic)=19$, the sextic~$\sextic$ is called \emph{maximizing}.
We emphasize that both
the inequality and the term apply to simple sextics only.

An irreducible sextic $\sextic\subset\Cp2$ is called \emph{special} (more precisely,
\emph{$\DG{2n}$-special}) if its fundamental group
$\pi_1:=\pi_1(\Cp2\sminus\sextic)$ factors to a dihedral group~$\DG{2n}$,
$n\ge3$.

A sextic $\sextic$ is
said to be
of \emph{torus type} if its
defining polynomial~$f$ can be written
in the form $f=f_2^3+f_3^2$,
where $f_2$ and $f_3$ are homogenous polynomials of degree 2 and 3, respectively.
A representation $f=f_2^3+f_3^2$ as above,
up to the obvious equivalence, is called a
\emph{torus structure} on~$\sextic$.
According to~\cite{degt:Oka},
an irreducible sextic $\sextic$
may have one,
four, or twelve distinct torus structures,
and we call~$\sextic$
a $1$-, $4$-, or $12$-torus sextic, respectively.
An irreducible sextic is of torus type if and only if it is $\DG6$-special,
see~\cite{degt:Oka}.
In this case, the group $\pi_1(\Cp2\sminus\sextic)$ factors to~$\MG$,
see~\cite{Zariski:group}.

The points of the intersection $f_2=f_3=0$ are singular for~$\sextic$; they are
called the \emph{inner} singularities of~$\sextic$ (with respect to the
given torus structure), whereas the other singular points are called
\emph{outer}. When listing the set of singularities of a $1$-torus sextic (or
describing a particular torus structure), it is common to enclose the inner
singularities in parentheses, \cf. \autoref{tab.torus}.
Conversely, the presence of a pair of parentheses in the notation indicates
that the sextic is of torus type.

Denote by $\all\cong\Cp{27}$ the space of all plane sextics. This space is
subdivided into equisingular strata $\all(\sset)$; we consider only those
with $\sset$ simple.
The space of all simple sextics and each of its
strata $\all(\sset)$ are further subdivided into families~$\all_*$,
$\all_*(\sset)$, where the subscript~$*$ refers to the sequence of invariant
factors of a certain finite group, see \autoref{s.htype} for the precise
definition.
Our primary concern are the spaces
\roster*
\item
$\all_\ns(\sset)$: non-special irreducible sextics, see \autoref{th.special},
and
\item
$\all_3(\sset)$: irreducible $1$-torus sextics, see \autoref{th.1-torus}.
\endroster
In this notation, irreducible $4$- and $12$-torus sextics constitute
$\all_{3,3}$ and $\all_{3,3,3}$, respectively, whereas irreducible
$\DG{2n}$-special sextics, $n=5,7$, constitute $\all_n$.
For each subscript~$*$, we denote by $\bar\all_*(\sset)$ and
$\partial\all_*(\sset):=\bar\all_*(\sset)\sminus\all_*(\sset)$
the closure and boundary of $\all_*(\sset)$
\emph{in~$\all_*$}.

\remark
The relation between torus type and the existence of certain dihedral
coverings (the families $\all_3$, $\all_{3,3}$, \etc.),
discovered for
irreducible sextics in~\cite{degt:Oka,degt:Oka2}
(see also Ishida, Tokunaga~\cite{Tokunaga:Oka} for reducible simple sextics),
is a manifestation of a much more general
phenomenon, \viz. a relation between the fundamental group of a
curve~$\sextic$ and
``special'' pencils containing~$\sextic$ (with an even further generalization
to quasi-projective varieties).
This was
studied in depth by E.~Artal, J.-I.~Cogolludo, A.~Libgober, and
others, see, \eg, recent
papers~\cite{Artal.Cogolludo.Libgober,Libgober.Cogolludo}.
\endremark

If $\sset$ is a simple set of singularities, the dimension of the
\emph{equisingular moduli space} $\all(\sset)/\!\PGL(3,\C)$ equals
$19-\Gm(\sset)$, as follows from the theory of
$K3$-surfaces.


The coordinatewise conjugation
$(z_0:z_1:z_2)\mapsto(\bar z_0:\bar z_1:\bar z_2)$ in $\Cp2$ induces a
real structure (\ie, anti-holomorphic involution)
$\conj\:\all\to\all$, which takes a sextic to its conjugate.
A sextic $\sextic\in\all$ is \emph{real} if $\conj(\sextic)=\sextic$. A
connected component $\Cal C\subset\all_*(\sset)$ is \emph{real} if
it is preserved by~$\conj$ as a set;
this property of~$\Cal C$ is independent of the choice of
coordinates in~$\Cp2$. Clearly, any
connected component containing a real curve is real.
The converse is not true; however, in the realm
of irreducible sextics, the only exception is $\all_\ns(\singset{A7+A6+A5})$,
see \autoref{prop.A7+A6+A5}.

Most results of the paper are stated
in terms of degenerations/perturbations of sets of
singularities and/or sextics (or, equivalently, in terms of adjacencies of
the equisingular strata of~$\all$).
As shown in~\cite{Looijenga:perturbations}, the deformation classes of
perturbations of a simple singular point~$P$ of type~$\sset$ are in a
one-to-one correspondence with
the isomorphism classes of primitive extensions $\sset'\pe\sset$ of root
lattices, see~\autoref{s.root} and~\autoref{s.extensions}.
Thus, by a \emph{degeneration} of sets of singularities we merely mean a
class of
primitive extensions $\sset'\pe\sset$ of root lattices. Recall
(see~\cite{Dynkin}) that $\sset'$ admits a degeneration to~$\sset$ if and
only if the Dynkin graph of~$\sset'$ is an induced subgraph of that
of~$\sset$.
A degeneration $\sextic'\dg\sextic$ of simple sextics gives rise to a
degeneration $\sset(\sextic')\pe\sset(\sextic)$.
According to~\cite{degt:8a2}, the converse
also holds: given a simple
sextic~$\sextic$ and a root lattice~$\sset'$, any degeneration
$\sset'\pe\sset(\sextic)$ is realized by a degeneration
$\sextic'\dg\sextic$ of simple sextics, so that $\sset(\sextic')=\sset'$.

\subsection{Lists and fundamental groups}\label{s.lists}
A complete list of the sets of singularities realized by simple plane sextics
is found in~\cite{Yang}, and the deformation classification of all maximizing
simple sextics is obtained in~\cite{Shimada:maximal}
(see also~\cite{degt:book} for an alternative approach to sextics with a
triple singular point).
The relevant part of these results is collected in Tables~\ref{tab.triple},
\ref{tab.double} (irreducible maximizing non-special sextics) and
\autoref{tab.torus} (irreducible maximizing \torus1 sextics).
In the tables, the column $(r,c)$ refers to the numbers of real ($r$) and
pairs of complex conjugate ($c$) curves realizing the given set of
singularities; thus, the total number of connected components of the stratum
$\all_\ns(\sset)$ (or $\all_3(\sset)$ for \autoref{tab.torus}) is $n:=r+2c$.
Some sets of singularities are prefixed with a link of the form
\tabrefform{n}: this link refers to the listings of the fundamental groups
found below.
Some pairs of singular points are marked with a $^*$. This marking is related
to the real structures; it is explained in \autoref{s.real.ns}.

%

\table
\caption{The spaces $\all_\ns(\sset)$, $\Gm(\sset)=19$,
with a triple point in $\sset$}\label{tab.triple}
\tabs\tab
2E8+A3//|(1,0)\cr
2E8+A2+A1//|(1,0)\cr
E8+E7+A4//|(0,1)\cr
E8+E7+2A2//|(1,0)\cr
E8+E6+D5//|(1,0)\cr
E8+E6+A5//|(0,1)\cr
E8+E6+A4+A1//|(1,0)\cr
E8+E6+A3+A2//|(1,0)\cr
E8+D11//|(1,0)\cr
E8+D9+A2//|(1,0)\cr
E8+D7+A4//|(1,0)\cr
E8+D5+A6//|(0,1)\cr
E8+D5+A4+A2//|(1,0)\cr
E8+A11//|(0,1)\cr
E8+A10+A1//|(1,1)\cr
E8+A9+A2//|(1,0)\cr
E8+A8+A3//|(1,0)\cr
E8+A8+A2+A1//|(1,1)\cr
E8+A7+A4//|(0,1)\cr
E8+A7+2A2//|(1,0)\cr
E8+A6+A5//|(0,1)\cr
E8+A6+A4+A1//|(1,1)\cr
E8+A6+A3+A2//|(1,0)\cr
E8+A6+2A2+A1//|(1,0)\cr
E8+A5+A4+A2//|(2,0)\cr
E8+A4+A3+2A2^*//pi1.E8+A4+A3+2A2|(1,0)\cr
E7+2E6^*//|(1,0)\cr
E7+E6+A6//|(0,1)\cr
E7+E6+A4+A2//|(2,0)\cr
E7+A12//|(0,1)\cr
E7+A10+A2//|(2,0)\cr
E7+A8+A4//|(0,1)\cr
E7+A6+A4+A2//|(2,0)\cr
E7+2A6//|(0,1)\cr
E7+2A4+2A2^*//pi1.E7+2A4+2A2|(1,0)\cr
2E6^*+A7//|(1,0)\cr
2E6^*+A6+A1//|(1,0)\cr
2E6+A4+A3//pi1.2E6+A4+A3|(1,0)\cr
E6+D13//|(1,0)\cr
E6+D11+A2//|(1,0)\cr
\tabbreak
E6+D9+A4//|(1,0)\cr
E6+D7+A6//|(1,0)\cr
E6+D5+A8//|(1,1)\cr
E6+D5+A6+A2//|(2,0)\cr
E6+D5+2A4//|(1,0)\cr
E6+A13//|(0,1)\cr
E6+A12+A1//|(0,1)\cr
E6+A10+A3//|(2,0)\cr
E6+A10+A2+A1//|(1,1)\cr
E6+A9+A4//|(1,1)\cr
E6+A8+A4+A1//|(1,1)\cr
E6+A7+A6//|(0,1)\cr
E6+A7+A4+A2//|(2,0)\cr
E6+A6+A4+A3//|(1,0)\cr
E6+A6+A4+A2+A1//|(1,1)\cr
E6+A5+2A4//|(2,0)\cr
D19//|(1,0)\cr
D17+A2//|(1,0)\cr
D15+A4//|(1,0)\cr
D13+A6//|(0,1)\cr
D13+A4+A2//|(1,0)\cr
D11+A8//|(1,0)\cr
D11+A6+A2//|(1,0)\cr
D11+A4+2A2^*//|(1,0)\cr
D9+A10//|(1,0)\cr
D9+A6+A4//|(1,0)\cr
D9+2A4^*+A2//|(1,0)\cr
D7+A12//|(1,1)\cr
D7+A10+A2//|(0,1)\cr
D7+A8+A4//|(2,0)\cr
D7+A6+A4+A2//|(1,0)\cr
D7+2A6//|(0,1)\cr
D5+A14//|(0,1)\cr
D5+A12+A2//|(1,0)\cr
D5+A10+A4//|(1,1)\cr
D5+A10+2A2^*//|(1,0)\cr
D5+A8+A6//|(0,1)\cr
D5+A8+A4+A2//|(1,1)\cr
D5+A6+2A4//|(2,0)\cr
D5+A6+A4+2A2^*//|(1,0)\cr
\endtab\endtabs

\endtable

The fundamental groups of most irreducible maximizing sextics are computed in
\cite{degt:book,degt:Artal};
the latest computations, using S.~Orevkov's recent
equations~\cite{Orevkov:equations}, are contained in~\cite{degt:equations}.
(Due to~\cite{Orevkov:equations}, the defining equations of
\emph{all} maximizing irreducible sextics with double points only are known
now.)
Quite a few sporadic computations of the fundamental groups are also
found
in~\cite{
Artal:trends,
Artal:invariants,
degt:e6,
degt:8a2,
degt:Oka3,
Oka:conjecture,
Oka:groups,
Oka.Pho:groups,
Zariski:9cuspidal} and a number of other papers, see~\cite{degt:book} for more
detailed references.

\table
\caption{The spaces $\all_\ns(\sset)$, $\Gm(\sset)=19$,
with double points only}\label{tab.double}
\tabs\tab
A19//|(2,0)\cr
A18+A1//|(1,1)\cr
A16+A3//|(2,0)\cr
A16+A2+A1//|(1,1)\cr
A15+A4//|(0,1)\cr
A14+A4+A1//pi1.unknown|(0,3)\cr
A13+A6//pi1.unknown|(0,2)\cr
A13+A4+A2//|(2,0)\cr
A12+A7//pi1.unknown|(0,1)\cr
A12+A6+A1//pi1.real|(1,1)\cr
A12+A4+A3//|(1,0)\cr
A12+A4+A2+A1//pi1.real|(1,1)\cr
A11+2A4^*//pi1.one|(2,0)\cr
A10+A9//|(2,0)\cr
A10+A8+A1//pi1.real|(1,1)\cr
\tabbreak
A10+A7+A2//|(2,0)\cr
A10+A6+A3//pi1.unknown|(0,1)\cr
A10+A6+A2+A1//pi1.real|(1,1)\cr
A10+A5+A4//|(2,0)\cr
A10+2A4^*+A1//pi1.unknown|(1,1)\cr
A10+A4+A3+A2//|(1,0)\cr
A10+A4+2A2+A1//|(2,0)\cr
A9+A6+A4//pi1.real|(1,1)\cr
A8+A7+A4//pi1.unknown|(0,1)\cr
A8+A6+A4+A1//pi1.real|(1,1)\cr
A7+2A6//pi1.unknown|(0,1)\cr
A7+A6+A4+A2//|(2,0)\cr
A7+2A4+2A2^*//|(1,0)\cr
2A6^*+A4+A2+A1//|(2,0)\cr
A6+A5+2A4^*//|(2,0)\cr
\endtab\endtabs

\endtable

The known fundamental groups $\pi_1:=\pi_1(\Cp2\sminus\sextic)$
of the maximizing non-special irreducible sextics~$\sextic$
are as follows (depending on the set of singularities):
\roster
\item\label{pi1.E8+A4+A3+2A2}
for
\singset{E8+A4+A3+2A2}, the group is
the central product
\[*
\pi_1=\SL(2,\FF5)\odot\CG{12}:=\bigl(\SL(2,\FF5)\times\CG{12}\bigr)/(-{\id}=6),
\]
where $-{\id}$
is the generator of the center $\CG2\subset\SL(2,\FF5)$;
\item\label{pi1.E7+2A4+2A2}
for
\singset{E7+2A4+2A2}, the group is
$\pi_1=\SL(2,\FF{19})\times\CG{6}$;
\item\label{pi1.2E6+A4+A3}
for
\singset{2E6+A4+A3}, the group is
$\pi_1=\SL(2,\FF5)\rtimes\CG6$, the generator of $\CG6$ acting on
$\SL(2,\FF5)$ by (any) order~$2$ outer automorphism;
\item\label{pi1.real}
for the six sets of singularities marked with \maketabref{pi1.real}
in \autoref{tab.double},
one has $(r,c)=(1,1)$, and \emph{only for the real curve}
the
group $\pi_1=\CG6$ is known;
\item\label{pi1.one}
for
\singset{A11+2A4},
only for one of the two curves the group $\pi_1=\CG6$ is known;
\item\label{pi1.unknown}
for the seven sets of singularities marked with \maketabref{pi1.unknown}
in \autoref{tab.double},
the fundamental group is still unknown.
\endroster
In all other cases, the fundamental group is abelian: $\pi_1=\CG6$.

The fundamental groups of sextics of torus type are large and more difficult
to describe. To simplify the description, we introduce a few \latin{ad hoc}
groups:
\[
G(\bar s):=\bigl\<\Ga_1,\Ga_2,\Ga_3\bigm|\Gr_3^4=(\Ga_1\Ga_2)^3,\
 \{\Ga_2\ra\Gs_1^i,\Ga_3\}_{s_i}=1,\ i=0,\ldots,5\bigr\>,
\label{gr.G}
\]
where $\bar s=(s_0,\ldots,s_5)\in\Z^6$ is an integral vector,
\begin{multline}
L_{p,q,r}:=\bigl\<\Ga_1,\Ga_2\bigm|(\Ga_1\Ga_2\Ga_1)^3=\Ga_2\Ga_1\Ga_2,\
 \{\Ga_2,(\Ga_1\Ga_2)\Ga_1(\Ga_1\Ga_2)\1\}_p\\
 =\{\Ga_1,\Ga_2\Ga_1\Ga_2\1\}_q
 =\{\Ga_2,(\Ga_1\Ga_2^2)\Ga_1(\Ga_1\Ga_2^2)\1\}_r=1\bigr\>,
\label{gr.L}
\end{multline}
where $p,q,r\in\Z$, and
\begin{multline}
E_{p,q}:=\bigl\<\Ga_1,\Ga_2,\Ga_3\bigm|
 \Gr_3\Ga_2\Gr_3\1=\Ga_2\1\Ga_1\Ga_2=\Gr_3\1\Ga_3\Gr_3,\\
 \Gr_3^4=(\Ga_1\Ga_2)^3,\
 \{\Ga_2,\Ga_3\}_p=\{\Ga_1,\Ga_3\}_q=1\bigr\>,
\label{gr.E}
\end{multline}
where $p,q\in\Z$.
Then,
the fundamental groups of the maximizing irreducible \torus1 sextics
are as
follows:
\roster
\item\label{pi1.(3E6)+A1}
for
\singset{(3E6)+A1},
the group is $\pi_1=\BG4/\Gs_2\Gs_1^2\Gs_2\Gs_3^3$;
\item\label{pi1.(2E6+A5)+A2}
for
\singset{(2E6+A5)+A2},
the groups are $E_{3,6}$, see~\eqref{gr.E},
and $L_{3,6,0}$, see~\eqref{gr.L};
\item\label{pi1.(2E6+2A2)+A3}
for
\singset{(2E6+2A2)+A3},
the group is $E_{4,3}$, see~\eqref{gr.E};
\item\label{pi1.(E6+A5+2A2)+A4}
for
\singset{(E6+A5+2A2)+A4},
the groups are $L_{5,6,3}$ and $G(6,5,3,3,5,6)$, see \eqref{gr.L}
and~\eqref{gr.G}, respectively;
\item\label{pi1.(A8+3A2)+A4+A1}
for
\singset{(A8+3A2)+A4+A1},
the group is
\begin{multline*}
\pi_1=\bigl\<\Ga_1,\Ga_2,\Ga_3\bigm|
[\Ga_2,\Ga_3]=\{\Ga_1,\Ga_2\}_3=\{\Ga_1,\Ga_3\}_9=1,\\
\Ga_3\Ga_1\Ga_2\1\Ga_3\Ga_1\Ga_3(\Ga_3\Ga_1)^{-2}\Ga_2=
 (\Ga_1\Ga_3)^2\Ga_2\1\Ga_1\Ga_3\Ga_2\Ga_1\bigr\>;
\end{multline*}
\item\label{pi1.(A8+A5+A2)+A4}
for the set of singularities \singset{(A8+A5+A2)+A4},
the group is unknown.
\endroster
In all other cases, the fundamental group is $\pi_1=\MG$.
In each of items~\ref{pi1.(2E6+A5)+A2} and~\ref{pi1.(E6+A5+2A2)+A4},
it is not known
whether the two groups are isomorphic.
The groups corresponding to distinct sets of singularities (listed above) are
distinct, except that it is not known whether the group in
\autoref{pi1.(A8+3A2)+A4+A1} is isomorphic to~$\MG$.

\table
\caption{The spaces $\all_3(\sset)$, $\Gm(\sset)=19$}\label{tab.torus}
\tabs\tab
(3E6)+A1//pi1.(3E6)+A1|(1,0)\cr
(2E6+A5)+A2//pi1.(2E6+A5)+A2|(2,0)\cr
(2E6+2A2^*)+A3//pi1.(2E6+2A2)+A3|(1,0)\cr
(E6+A11)+A2//|(1,0)\cr
(E6+A8+A2)+A3//|(1,0)\cr
(E6+A8+A2)+A2+A1//|(1,1)\cr
(E6+A5+2A2^*)+A4//pi1.(E6+A5+2A2)+A4|(2,0)\cr
D5+(A8+3A2^*)//|(1,0)\cr
\tabbreak
(A17)+A2//|(1,0)\cr
(A14+A2)+A3//|(1,0)\cr
(A14+A2)+A2+A1//|(1,0)\cr
(A11+2A2^*)+A4//|(1,0)\cr
(2A8)+A3//|(1,0)\cr
(A8+A5+A2)+A4//pi1.(A8+A5+A2)+A4|(0,1)\cr
(A8+3A2^*)+A4+A1//pi1.(A8+3A2)+A4+A1|(1,0)\cr
\omit\strut\cr
\endtab\endtabs

\endtable

\subsection{Statements}\label{s.statements}
There are
$110$
maximizing sets of simple singularities
realized
by non-special irreducible sextics.
We found that $2996$ sets of simple singularities are realized by
non-maximizing non-special irreducible sextics.
(This statement is almost contained in~\cite{Yang}, although no distinction
between special and non-special curves is made there,
nor a description of \emph{non-maximizing} irreducible sextics.)
The corresponding counts for
irreducible \torus1 sextics are $15$ and $105$, respectively,
see~\cite{Oka.Pho:moduli}. Our principal results (the deformation
classification and a few consequences on the fundamental group) are stated in
the rest of this section, with references to the proofs given in the headers.


\table
\caption{Disconnected spaces $\all_\ns(\sset)$, $\Gm(\sset)<19$}\label{tab.disconnected}
\tabs\tab
E8+2A5//|(0,1)\cr
E7+E6+A5//|(0,1)\cr
E7+A7+A4//|(0,1)\cr
E6+A11+A1//|(0,1)\cr
E6+A7+A5//|(0,1)\cr
E6+A6+A5+A1//|(0,1)\cr
E6+2A5+A1//|(0,1)\cr
\tabbreak
D6+2A6//|(0,1)\cr
D5+2A6+A1//|(0,1)\cr
2A9//|(2,0)\cr
A7+A6+A5//|(1,0)\cr
3A6//|(0,1)\cr
2A6+2A3//|(0,1)\cr
2A7+A4//|(0,1)\cr
\endtab\endtabs

\endtable

\theorem[see \autoref{proof.classification} and \autoref{s.real}]\label{th.classification}
The space $\all_\ns(\sset)$ is nonempty
if and only if
either $\sset$ is in one of the following two exceptional degeneration chains
\[*
\singset{2D8}\dg\singset{D9+D8}\dg\singset{2D9},\qquad
\singset{2D4+4A2}\dg\singset{D7+D4+3A2}\dg\singset{2D7+2A2}
\]
or $\sset$ degenerates to one of the maximizing
sets
of singularities
listed in Tables~\ref{tab.triple},~\ref{tab.double}.
The numbers $(r,c)$ of connected components of $\all_\ns(\sset)$
are as shown in
Tables~\ref{tab.triple}, \ref{tab.double}, and~\ref{tab.disconnected}\rom;
in all other cases, $\all_\ns(\sset)$ is connected and contains real curves.
\endtheorem

Two lines in \autoref{tab.disconnected} deserve separate statements:
to our knowledge, phenomena of this kind have not been observed before.

\proposition[see \autoref{proof.2A9}]\label{prop.2A9}
Let
$\sset_0:=\singset{2A9}$, $\sset_1:=\singset{A19}$, and
$\sset_2:=\singset{A10+A9}$. The space $\all_\ns(\sset_i)$, $i=0,1,2$,
consists of two connected components $\all^\pm_\ns(\sset_i)$,
each containing real curves,
so that
$\partial\all^\epsilon_\ns(\sset_0)=
 \all^\epsilon_\ns(\sset_1)\cup\all^\epsilon_\ns(\sset_2)$ for
each $\epsilon=\pm$.
\endproposition

\proposition[see \autoref{proof.A7+A6+A5}]\label{prop.A7+A6+A5}
The space $\all_\ns(\singset{A7+A6+A5})=\all(\singset{A7+A6+A5})$ is
connected \rom(hence, its only component is real\rom),
but it contains no real curves.
\endproposition

In the other cases in \autoref{tab.disconnected},
the space
$\all_\ns(\sset)$ consists of two
complex
conjugate components. The first such example, \viz.
$\sset=\singset{E6+A11+A1}$, was discovered in~\cite{Aysegul:paper}.
The adjacencies of these non-real components are studied in
\autoref{s.adjacencies}. Note that one set of singularities, \viz.
\singset{E6+2A5+A1}, has Milnor number~$17$; it gives rise to an interesting
adjacency phenomenon, see \autoref{prop.p=3}.

\corollary[see \autoref{proof.degeneration}]\label{cor.degeneration}
With the same six exceptions as in \autoref{th.classification},
any non-special irreducible simple sextic degenerates to a
maximizing sextic with these properties,
see Tables~\ref{tab.triple} and~\ref{tab.double}.
\endcorollary

\corollary[see \autoref{proof.pi1}]\label{cor.pi1}
Let~$\sextic\subset\Cp2$ be a non-special irreducible simple plane sextic.
If $\Gm(\sextic)=19$, the fundamental group
$\pi_1:=\pi_1(\Cp2\sminus\sextic)$ is as shown in
Tables~\ref{tab.triple} and~\ref{tab.double}. Otherwise, one has
\roster*
\item
$\pi_1=\SL(2,\FF3)\times\CG2$ for
\singset{2D7+2A2}, \singset{D7+D4+3A2}, and \singset{2D4+4A2},
\item
$\pi_1=\SL(2,\FF5)\odot\CG{12}$, see \autoref{s.lists},
for \singset{2A4+2A3+2A2},
\endroster
and $\pi_1=\CG6$ in all other cases.
\endcorollary

The remaining statements deal with sextics of torus type, and we introduce
the notion of weight.
The \emph{weight} $w(P)$ of a simple singular point~$P$
is defined
\via\
$w(\bA_{3p-1})=p$, $w(\bE_6)=2$, and $w(P)=0$ otherwise. The weight of a set
of singularities~$\sset$ (or a simple sextic~$\sextic$) is the total weight
of its singular points. Recall
(see~\cite{degt:Oka}) that, if $\sextic$ is a \torus1 sextic, then
$6\le w(\sextic)\le7$.
Conversely, if $\sextic$ is an irreducible sextic and
either $w(\sextic)=7$ or $w(\sextic)=6$
and $\sextic$ has at least one singular point
$P\ne\bA_1$ of weight~$0$, then $\sextic$ is a \torus1 sextic.

\theorem[see \autoref{proof.torus} and \autoref{s.real.torus}]\label{th.torus}
A set of singularities~$\sset$ with $w(\sset)\ge6$
is realized by an irreducible simple
\torus1 sextic~$\sextic$
if and only if $\sset$ degenerates to one of the maximizing
sets listed in \autoref{tab.torus}.
Furthermore,
if $\Gm(\sset)\le18$, a sextic~$\sextic$ as above is unique up to
equisingular deformation
and the space $\all_3(\sset)$ contains real curves.
\endtheorem

\corollary[see \autoref{proof.torus}]\label{cor.torus}
Any irreducible simple \torus1 sextic
degenerates to a
maximizing sextic with these properties,
see \autoref{tab.torus}.
\endcorollary

There are $51$ sets of singularities~$\sset$ (all of
weight~$6$) realized by both \torus1 and non-special irreducible sextics.
Formally, these sets of singularities can be extracted from
Theorems~\ref{th.classification} and~\ref{th.torus}; an explicit list is found
in~\cite{Aysegul:paper}.
The corresponding sextics constitute the so-called \emph{classical Zariski
pairs}.

\corollary[see \autoref{proof.pi1.torus}]\label{cor.pi1.torus}
Let $\sextic\subset\Cp2$ be an irreducible simple \torus1 sextic.
If $\Gm(\sextic)=19$, the fundamental group
$\pi_1:=\pi_1(\Cp2\sminus\sextic)$ is as shown in \autoref{tab.torus}.
Otherwise, one has
$\pi_1=\BG4/\Gs_2\Gs_1^2\Gs_2\Gs_3^3$ for
the sets of singularities
\[*
\gathered
\singset{(2E6+2A2)+2A1},\quad
\singset{(E6+4A2)+3A1},\quad
\singset{(E6+4A2)+A3+A1},\\
\singset{(6A2)+A3+2A1},\quad
\singset{(6A2)+4A1},
\endgathered
\]
and $\pi_1=\MG$ in all other cases.
\endcorollary

\remark
In \autoref{cor.pi1.torus}, the non-maximizing \torus1 sextics with the group
$\pi_1=\BG4/\Gs_2\Gs_1^2\Gs_2\Gs_3^3$ can be characterized as the
degenerations of \singset{(6A2)+4A1}.
\endremark

\subsection{Other irreducible sextics}\label{s.special}
For the reader's convenience and completeness of the exposition, we recall
the classification
of
the other irreducible simple sextics, \viz.
the $\DG{10}$- and $\DG{14}$-special sextics and the $4$- and $12$-torus
ones.
The fundamental groups are computed in several papers, see~\cite{degt:book}
for detailed references.

\theorem[see~\cite{degt:Oka}]\label{th.5-special}
The space $\all_5$ consists of eight connected components,
one component $\all_5(\sset)$ for
each of the following sets of singularities~$\sset$\rom:
\[*
\gathered
\singset{2A9},\quad
\singset{A9+2A4+A2},\quad
\singset{A9+2A4+A1},\quad
\singset{A9+2A4},\\
\singset{4A4+A2},\quad
\singset{4A4+2A1},\quad
\singset{4A4+A1},\quad
\singset{4A4}.
\endgathered
\]
All components are real and contain real curves.
\pni
\endtheorem

The fundamental group $\pi_1:=\pi_1(\Cp2\sminus\sextic)$ of a
\emph{simple} sextic
$\sextic\in\all_5(\sset)$ can be described as follows.
Denoting temporarily by~$G'$ the
derived subgroup $[G,G]$,
one always has $\pi_1/\pi_1''=\DG{10}\times\CG3$. Besides,
\roster
\item\label{i.5.1}
if $\sset=\singset{A9+2A4+A2}$, then
$\pi_1''$ is the only perfect group of order~$120$;
\item\label{i.5.2}
if $\sset=\singset{4A4+2A1}$, then
$\pi_1''/\pi_1'''=\CG2^4$ and $\pi_1'''=\CG2$\rom, so that
$\ord\pi_1=960$\rom;
\item
in all other cases, $\pi_1=\DG{10}\times\CG3$.
\endroster
The precise presentations
in~\iref{i.5.1} and~\iref{i.5.2}
are rather lengthy, and we refer to~\cite{degt:Oka3}.

\theorem[see~\cite{degt:Oka}]\label{th.7-special}
The space $\all_7$ consists of two connected components,
one component $\all_7(\sset)$ for
each of the following sets of singularities~$\sset$\rom:
\[*
\singset{3A6+A1},\quad\singset{3A6}.
\]
Both components are real and contain real curves.
\pni
\endtheorem

The fundamental groups of all $\DG{14}$-special sextics are
$\DG{14}\times\CG3$.

\remark
The sets of singularities
\singset{2A9},
\singset{A9+2A4+A1},
\singset{A9+2A4},
\singset{4A4+A1},
\singset{4A4} (\cf. \autoref{th.5-special})
and \singset{3A6} (\cf. \autoref{th.7-special})
are also realized by non-special irreducible sextics,
each by a single connected deformation family.
\endremark

\theorem[see~\cite{degt:Oka}]\label{th.>=8}
The union $\all_{3.3}\cup\all_{3,3,3}$ consists of eight connected components,
one component for
each of the following sets of singularities~$\sset$\rom:
\roster*
\item
$\all_{3,3}$ \rom($4$-torus sextics,
\latin{idem} weight~$w=8$\rom)\rom:
\singset{E6+A5+4A2},
\singset{E6+6A2},
\singset{2A5+4A2},
\singset{A5+6A2+A1},
\singset{A5+6A2},
\singset{8A2+A1},
\singset{8A2}\rom;
\item
$\all_{3,3,3}$ \rom($12$-torus sextics, \latin{idem} $w=9$\rom)\rom:
\singset{9A2}.
\endroster
All components are real and contain real curves.
\pni
\endtheorem

All sets of singularities of weight~$8$ degenerate to \singset{E6+A5+4A2} and
can be characterized as perturbations of the latter preserving all four torus
structures. Note that \singset{9A2}
does \emph{not} degenerate to a maximizing sextic, irreducible or
not!

Introduce the group
\[*
\HH:=\bigl\<\Ga,\bar\Ga,\Gb,\Gg,\bar\Gg\bigm|
 \{\Ga,\Gb\}_3=\{\bar\Ga,\Gb\}_3=
 \{\Gg,\Gb\}_3=\{\bar\Gg,\Gb\}_3=
 \Gb\Gg\Ga\Gb\bar\Gg\bar\Ga=1\bigr\>.
\]
In this notation (see also~\eqref{gr.G}),
the fundamental group $\pi_1:=\pi_1(\Cp2\sminus\sextic)$ of a
sextic $\sextic$ with a set of singularities~$\sset$ of weight~$8$ or~$9$
is as follows:
\roster
\item
if $\sset=\singset{9A2}$ ($w=9$), then
\[*
\pi_1=H_3:=\HH/\quo
 {\{\bar\Gg,\Ga\}=\{\Gg,\bar\Ga\}=\{\Gg,\bar\Gg\}=
 [\Gb,\Ga\1\Gg\1\bar\Ga\bar\Gg]=1,\
 \bar\Gg\1\Ga\bar\Gg=\Gg\1\bar\Ga\Gg};
\]
\item
if $\sset=\singset{E6+A5+4A2}$, then
\[*
\pi_1=H_2:=\HH/
 \quo{\bar\Ga\Gg\Ga=\Ga\bar\Gg\bar\Ga=\Gg\Ga\bar\Gg=\bar\Gg\bar\Ga\Gg}
 \cong G(3,6,3,3,6,3);
\]
\item
if $\sset=\singset{A5+6A2+A1}$, then
\[*
\pi_1=H_1:=\HH/\quo{
 \{\Ga,\Gg\}_3=\{\bar\Ga,\bar\Gg\}_3=[\Gg,\bar\Gg]=1,\
 \Gg\Ga\bar\Gg=\bar\Gg\bar\Ga\Gg};
\]
\item
for all other sextics of weight~$8$,
\[*
\pi_1=H_0:=\HH/\quo{\Ga=\bar\Ga,\Gg=\bar\Gg,\{\Ga,\Gg\}_3=1}
 \cong G(3,3,3,3,3,3).
\]
\endroster
All perturbation epimorphisms
$H_3\onto H_0$ and $H_2\onto H_1\onto H_0$, \cf. \autoref{th.Zariski},
lift to the identity $\HH\to \HH$. We do not know whether
the epimorphism
$H_2\onto H_1$ is proper; the others are.

\subsection{Further generalizations}\label{s.reducible}
Altogether, there are $11308$ configurations (in the sense
of~\cite{degt:JAG}) of simple sextics, irreducible or not.
This result was first announced in~\cite{Shimada:splitting}, where
configurations are called \emph{lattice types}; roughly, these are
certain sets of
lattice
data
invariant under equisingular deformations and
recording both the position of the
singularities with respect to the irreducible components of the curve and the
existence of dihedral coverings.

The corresponding equisingular strata split into $11272$ real and $132$ pairs
of complex conjugate components.
As expected, this discrepancy is mainly due to the maximizing curves
(\latin{ergo} definite transcendental lattices), see~\cite{Shimada:maximal};
if $\mu<19$, then, in addition to \autoref{tab.disconnected}, there is a
single stratum $\all_2(\singset{2A9})$ consisting of two real
components (the sextic splits into an irreducible quintic and a line) and ten strata
(eight sets of singularities)
consisting of pairs of complex conjugate components.
Furthermore,
the stratum $\all_\ns(\singset{A7+A6+A5})$ remains the only real connected
component not containing real curves, \cf.~\autoref{prop.A7+A6+A5}.

There are $629$ maximizing configurations
($\mu=19$; see \cite{Shimada:maximal,Yang}). Besides,
there are
$16$ (with $\mu=18$) and $2$ (with $\mu=17$) other configurations
extremal with respect to degeneration
(\cf. the existence part of \autoref{th.classification}).
A thorough analysis of the adjacencies of the strata
and the computation of various symmetry groups
(in particular,
analogues of \autoref{th.group} and \autoref{s.adjacencies}
for reducible curves)
still require some work; therefore, we postpone the details until a later
paper.

\section{Integral lattices}\label{S.lattices}

\subsection{Finite quadratic forms\noaux{ (see~\cite{Nikulin:forms})}}\label{s.forms}
A \emph{finite quadratic form} is a finite abelian group $\dN$ equipped with a
symmetric bilinear form $b\:\dN\otimes\dN\to\Q/\Z$ and a
\emph{quadratic extension}
of~$b$, \ie, a map $q\:\dN\to\Q/2\Z$ such that
$q(x+y)-q(x)-q(y)=2b(x,y)$ for all $x,y \in \dN$
(where $2$ is the isomorphism ${\times2}\:\Q/\Z\rightarrow \Q/2\Z$);
clearly, $b$ is determined by~$q$.
If $q$ is understood, we abbreviate
$b(x,y)=x\cdot y$ and $q(x)=x^2$.
In what follows, we consider \emph{nondegenerate}
forms only, \ie, such that the associated homomorphism
$\dN\to\Hom(\dN,\Q/\Z)$, $x\mapsto(y\mapsto x\cdot y)$ is an isomorphism.

Each
finite quadratic
form $\dN$ splits into orthogonal sum
$\dN=\bigoplus_{p\in\PP} \dN_p$
of its $p$-primary
components
$\dN_p:=\dN\otimes\Z_p$.
The \emph{length} $\ell(\dN)$ of~$\dN$
is the minimal number of generators of $\dN$.
Obviously, $\ell(\dN)=\max_{p\in\PP}\ell_p(\dN)$, where
$\ell_p(\dN):=\ell(\dN_p)$.
The notation $-\dN$ stands for the group~$\dN$ with the form $x\mapsto-x^2$.

We describe
nondegenerate
finite quadratic forms by expressions of the form
$\<q_1\>\oplus\ldots\oplus\<q_r\>$, where $q_i:=\frac{m_i}{n_i}\in\Q$,
$\gcd(m_i,n_i)=1$, $m_in_i=0\bmod2$; the group is generated by pairwise
orthogonal elements $\Ga_1,\ldots,\Ga_r$ (numbered in the order of
appearance), so that $\Ga_i^2=q_i\bmod2\Z$ and the order of~$\Ga_i$
is~$n_i$. (In the $2$-torsion, there also may be indecomposable summands of
length~$2$, but we do not need them.)
Describing an automorphism~$\Gs$ of such a group,
we only list the images of the generators~$\Ga_i$ that are moved by~$\Gs$.

A finite quadratic form is called \emph{even} if $x^2=0\bmod\Z$ for each
element
$x\in\dN$ of order two; otherwise, the form is called \emph{odd}.
In other words, $\dN$ is odd if and only it contains $\<\pm\frac{1}{2}\>$
as an orthogonal summand.

Given a prime $p\in\PP$,
the \emph{determinant} $\det_p\dN$
is defined as the determinant of the
`matrix' of the quadratic form on~$\dN_p$ in an appropriate basis
(see~\cite{Miranda.Morrison:book} for the technical details);
\eg,
it is sufficient, although not necessary,
to take for a basis the set of generators of the
indecomposable cyclic (and those of length~$2$ if $p=2$) summands
constituting an orthogonal decomposition of~$\dN_p$.
Alternatively, $\det_p\dN$ is originally defined in~\cite{Nikulin:forms}
as the determinant of the unique $p$-adic lattice~$N_p$
such that $\rank N_p=\ell(\dN_p)$ and $\discr N_p=\dN_p$.
The determinant is an element of $\Q_p$ well defined modulo
the group of squares $(\Q_p\units)^2$; if $\dN_p$ is nondegenerate, then
one has
$\det_p\dN=u/\ls|\dN_p|$ for some $u\in\Z_p\units/(\Z_p\units)^2$.
In the case $p=2$, the determinant $\det_2\dN$ is
well defined only if $\dN_2$ is even (as otherwise a $2$-adic lattice $N_2$
as  above is not unique: there are two isomorphism classes whose determinants
differ by $5\in\Z_2\units$).
By definition,
one always has $\ls|\dN|\det_p\dN\in\Z_p\units/(\Z_p\units)^2$.

The group of $q$-autoisometries
of~$\dN$ is denoted by $\Aut\dN$;
obviously, one has $\Aut\dN=\prod_{p\in\PP}\Aut\dN_p$.
An element $\xi\in\dN_p$ is called a \emph{mirror}
if, for some integer~$k$, one has $p^k\xi=0$ and
$\xi^2=2u/p^k\bmod2\Z$, $\gcd(u,p)=1$.
If this is the case, the map $x\mapsto2(x\cdot\xi)/\xi^2\bmod p^k$
is a well defined
functional $\dN_p\to\Z/p^k$; hence, one has a
\emph{reflection} $\refl\xi\in\Aut\dN_p$,
\[*
\refl\xi\:x\mapsto x-\frac{2(x\cdot\xi)}{\xi^2}\xi.
\]
Note that $\refl\xi=\id$ whenever
$2\xi=0$ and $\xi^2=\frac12\bmod\Z$.

\subsection{Lattices and discriminant forms\noaux{ (see~\cite{Nikulin:forms})}}\label{s.lattices}
An \emph{\rom(integral\rom) lattice} $N$
is a finitely generated free abelian group equipped with a symmetric
bilinear form $b\: N\otimes N \rightarrow \mathbb{Z}$.
If $b$ is understood, we abbreviate
$b(x,y)=x\cdot y$ and
$b(x,x)=x^{2}$. A lattice $N$ is called
\emph{even} if $x^{2}=0 \bmod 2$ for all $x \in N$;
it is called
\emph{odd} otherwise.
The determinant $\det N$ of a lattice $N$ is
the
determinant of the Gram matrix
of $b$. As the transition matrix from one integral basis to another
has
determinant $\pm 1$,
the determinant
$\det N \in \mathbb{Z}$ is well-defined. The lattice $N$
is called \emph{non-degenerate} if $\det N \neq 0$ and \emph{unimodular}
if $\det N = \pm 1$. The signature $(\sigma_{+} N, \sigma_{-} N)$ of a
non-degenerate lattice $N$ is
the pair of
the inertia indices of the bilinear form $b$.

For a lattice $N$, the bilinear form extends to a $\mathbb{Q}$-valued bilinear
form on $N \otimes \Q$. If $N$ is non-degenerate, the dual group
$N\dual:=\Hom(N,\Z)$ can be identified with the subgroup
$\{x \in N\otimes\Q \mid \text{$x\cdot y \in \Z$ for all $y \in N$}\}$.
The lattice $N$ is a finite index subgroup of $N\dual$.
The quotient $\discr N:= N\dual/N$ is called the \emph{discriminant group}
of~$N$;
it is often
denoted by $\dN$, and we use the shortcut $\discr_pN=\dN_p$ for the
$p$-primary components.
One has
$\det N=(-1)^{\Gs_-N}\ls|\dN|$.
The group~$\dN$ inherits from $N\otimes\Q$ a symmetric bilinear
form $b\:\dN\otimes\dN\rightarrow\Q/\Z$, called
the \emph{discriminant form}, and, if $N$ is even,
a quadratic extension of $b$.

\convention
Unless
specified otherwise, all lattices considered below are nondegenerate
and even. The discriminant group of such a lattice is always
regarded as a finite quadratic form.
\endconvention

The
\emph{genus} $\genus(N)$ of a
nondegenerate even lattice~$N$ can be defined as the set of
isomorphism classes of all even lattices~$L$ such that
$\discr L\cong\dN$ and $\Gs_\pm L=\Gs_\pm N$.
If $N$ is indefinite and $\rank N\ge3$, then $\genus(N)$ is a
finite abelian group with
the group law given by \autoref{th.MM} below.

An \emph{isometry} of lattices
is a homomorphism of
abelian groups preserving the forms. (Note that we do not assume the
surjectivity.)
The group of auto-isometries of a lattice~$N$ is denoted by $\OG(N)$.
There is an obvious natural homomorphism $\map\:\!\OG(N)\to\Aut\dN$, and we denote by
$\map_p\:\!\OG(N)\to\Aut\dN_p$ its restrictions to the $p$-primary components.
For an element $u\in N$ such that $2u/u^2\in N\dual$,
the \emph{reflection} $\refl{u}\:x\mapsto2u(x\cdot u)/u^2$ is an involutive
isometry of~$N$. Each image $\map_p(\refl{u})$, $p\in\PP$, is also a
reflection. If $u^2=\pm1$ or~$\pm2$, then $\map(\refl{u})=\id$.


\subsection{Root lattices\noaux{ (see~\cite{Bourbaki:Lie})}}\label{s.root}
In this paper, a \emph{root lattice} is a negative definite lattice generated
by vectors of square $(-2)$ (\emph{roots}). Any root lattice has a unique
decomposition into orthogonal sum of indecomposable ones, which are of types
$\bA_p$, $p\ge1$, $\bD_q$, $q\ge4$, $\bE_6$, $\bE_7$, or~$\bE_8$.

Given a root lattice~$S$,
the vertices of the Dynkin graph $\Dynkin:=\Dynkin_S$
can be identified with the elements of a basis for~$S$
constituting a single Weyl chamber. This identification defines a
homomorphism $\Sym\Dynkin\to\OG(S)$, $s\mapsto s_*$, where $\Sym\Dynkin$ is
the group of symmetries of~$\Dynkin$. The image consists of the isometries
preserving the distinguished Weyl chamber.
For indecomposable root lattices, the groups $\Sym\Dynkin$ are as follows:
\roster*
\item
$\Sym\Dynkin=1$ if $S$ is $\bA_1$, $\bE_7$, or $\bE_8$,
\item
$\Sym\Dynkin\cong\SG3\cong\DG6$ if $S$ is $\bD_4$, and
\item
$\Sym\Dynkin=\CG2$ in all other cases.
\endroster
In the latter case, unless $S=\bD\subtext{even}$, the generator of
$\Sym\Dynkin$ induces $-{\id}$ on
the discriminant $\dS:=\discr S$. If $S=\bE_8$, then $\dS=0$.
For $S=\bA_1$, $\bE_7$, or~$\bD\subtext{even}$, the
discriminant
groups $\dS$ are
$\FF2$-modules and $-{\id}=\id$ in $\Aut\dS$.

A choice of a Weyl chamber gives rise to a decomposition
$\OG(S)=R(S)\rtimes\Sym\Dynkin$, where
$R(S)\subset\OG(S)$ is the subgroup generated by
reflections $\refl{u}$, $u\in S$, $u^2=-2$.
Furthermore,
\[*
\Ker[\map\:\!\OG(S)\to\Aut\dS]=R(S)\rtimes\Sym_0\Dynkin,
\]
where $\Sym_0\Dynkin$ is the group of permutations of the $\bE_8$-type
components of~$\Dynkin$.
Thus, denoting by $\Sym'\Dynkin\subset\Sym\Dynkin$ the group of
symmetries acting identically on the union of the $\bE_8$-type components,
we obtain an isomorphism $\Sym'\Dynkin=\Im\map$.
For future references, we combine these statements in a separate lemma.

\lemma\label{lem.root}
Let $S$ be a root lattice. Then, the epimorphism $\map\:\!\OG(S)\onto\Im\map$
has a splitting $\Im\map=\Sym'\Dynkin_S\into\OG(S)$, and one always has
$-{\id}\in\Im\map$.
\done
\endlemma


\subsection{Lattice extensions\noaux{ (see~\cite{Nikulin:forms})}}\label{s.extensions}
An
\emph{extension} of a lattice $S$ is
an isometry $S\to L$.
Two extensions
$S\to L_1,L_2$ are \emph{\rom(strictly\rom) isomorphic} if there is a
bijective isometry
$L_1\to L_2$ identical on~$S$.
More generally, given a subgroup $\OG'\subset\OG(S)$, two extensions are
\emph{$\OG'$-isomorphic} if they are related by a bijective isometry whose
restriction to~$S$ is an element of~$\OG'$.

We use the notation $S\fe L$ for finite index extension ($[L:S]<\infty$).
There is
a unique embedding $L\subset S\otimes\Q$ and, hence,
inclusions $S\subset L\subset L\dual\subset S\dual$.
The \emph{kernel} of a finite index extension $S\fe L$ is the subgroup
$\CK:=L/S\subset S\dual/S=\dS$. Since $L$ is an even integral lattice,
the kernel $\CK$ is isotropic, \ie, the restriction to $\CK$ of
the quadratic form $q:\dS\to\Q/2\Z$ is
identically zero. Conversely, given an isotropic subgroup $\CK\subset\dS$,
the subgroup $L=\{ x\in S\dual\,|\,(x \bmod S)\in\CK\}\subset S\dual$
is an extension of $S$. Thus, we have the following theorem.

\theorem[Nikulin~\cite{Nikulin:forms}]\label{th.N.fin.ind.extensions}
The map $L\mapsto\CK=L/S\subset\dS$ establishes a one-to-one correspondence
between the set of
strict
isomorphism classes of finite index extension $S\fe L$
and that of isotropic subgroup $\CK\subset\dS$. One has $\dL=\CK^\perp\!/\CK$.
\pni
\endtheorem

An isometry $a\in\OG(S)$ extends to a finite index extension~$L$ if and only
if $\map(a)$ preserves the kernel~$\CK$ (as a set).
Hence, $\OG'$-isomorphism classes of finite index extensions of~$S$
correspond to the $\map(\OG')$-orbits of isotropic subgroups $\CK\subset\dS$.

Another extreme case is that of a \emph{primitive} extension $S\to L$, \ie,
such that the group $L/S$ is torsion free; we use the notation $S\pe L$.
If $L$ is unimodular, one has $\discr S^\perp\cong-\dS$: the graph of this
anti-isometry is the kernel of the finite index extension
$S\oplus S^\perp\fe L$.
Hence,
the genus $\genus(S^\perp)$ is
determined by those of~$S$ and~$L$.
If $L$ is also indefinite, it is unique in its genus. Then,
for each
representative $N\in\genus(S^\perp)$,
an extension $S\pe L$ with $S^\perp\cong N$ is determined by a
bijective anti-isometry $\Gf\:\dS\to\dN$
($L$ is the finite index extension of $S\oplus N$ whose kernel is
the graph of~$\Gf$),
and the latter induces a homomorphism
$\map^\Gf\:\!\OG(S)\to\Aut\dN$. If $\Gf$ is not fixed, this map is well
defined up to an inner automorphism of $\Aut\dN$.

\theorem[Nikulin~\cite{Nikulin:forms}]\label{th.N.extensions}
Let~$L$ be an indefinite unimodular even lattice, $S\subset L$ a
nondegenerate primitive sublattice, and $\OG'\subset\OG(S)$ a subgroup.
Then, the $\OG'$-isomorphism classes of primitive extensions $S\pe L$
are enumerated by
the pairs $(N,c_N)$, where $N\in\genus(S^\perp)$ and
$c_N\in\map^\Gf(\OG')\backslash\Aut\dN/\Im\map$ is a double coset
\rom(for given~$N$ and some anti-isometry $\Gf\:\dS\to\dN$\rom).
\pni
\endtheorem

\theorem[Nikulin~\cite{Nikulin:forms}]\label{th.N.isometry}
Let $S\pe L$ be a lattice extension as in \autoref{th.N.extensions},
$N=S^\perp$, and $\Gf\:\dS\to\dN$ the corresponding anti-isometry. Then, a
pair of isometries $a_S\in\OG(S)$, $a_N\in\OG(N)$ extends to~$L$ if and only
if $\map^\Gf(a_S)=\map(a_N)$.
\pni
\endtheorem

Fix the notation $\bL:=2\bE_8\oplus3\bU$, where $\bU$ is the \emph{hyperbolic
plane}, $\bU=\Z u+\Z v$, $u^2=v^2=0$, $u\cdot v=1$,
and $\bE_8$ is the root lattice, see \autoref{s.root}.
For the ease of references, we recast Nikulin's existence theorem
from~\cite{Nikulin:forms} to the particular case of primitive extensions
$S\pe\bL$. Note that we do not need the restriction on the Brown invariant:
by the additivity, it would hold automatically.

\theorem[Nikulin~\cite{Nikulin:forms}]\label{th.N.existence}
Given a nondegenerate even lattice~$S$, a primitive extension $S\pe\bL$ exists
if and only if
the following conditions hold\rom:
$\Gs_+S\le3$, $\Gs_-S\le19$, $\ell(\dS)\le\Gd:=22-\rank S$, and
\roster*
\item
for all odd $p\in\PP$, either $\ell_p(\dS)<\Gd$ or
$\ls|\dS|\det_p\dS=(-1)^{\Gs_+S-1}\bmod(\Z_p\units)^2$\rom;
\item
either $\ell_2(\dS)<\Gd$, or $\dS_2$ is odd, or
$\ls|\dS|\det_2\dS=\pm1\bmod(\Z_2\units)^2$.
\pni
\endroster
\endtheorem

\subsection{Miranda--Morrison results\noaux{
(see~\cite{Miranda.Morrison:1,Miranda.Morrison:2,Miranda.Morrison:book})}}\label{s.MM}
%
Classically,
the uniqueness of a lattice~$N$ in its genus and the surjectivity of the map
$\map\:\OG(N)\to\Aut\dN$ are established using the sufficient conditions
found in~\cite{Nikulin:forms}. Unfortunately, these results do not cover our
needs, and we use the stronger criteria developed
in~\cite{Miranda.Morrison:1,Miranda.Morrison:2,Miranda.Morrison:book}.
Throughout the rest of this section, we assume that
\roster*
\item[$(*)$]\hypertarget{convention}
$N$ is a nondegenerate indefinite even lattice, $\rank N\ge3$.
\endroster

\warning
The convention used in this paper (following Nikulin~\cite{Nikulin:forms}
and, eventually, Gauss) differs slightly from that of Miranda--Morrison,
where quadratic and bilinear forms are related \via\
$q(x+y)-q(x)-q(y)=b(x,y)$.
Roughly, the values of all quadratic (but not bilinear) forms
in~\cite{Miranda.Morrison:1,Miranda.Morrison:2,Miranda.Morrison:book},
both on lattices and finite groups, should be multiplied by~$2$. In
particular, all lattices
in~\cite{Miranda.Morrison:1,Miranda.Morrison:2,Miranda.Morrison:book}
are even by definition.
Note though that this multiplication by~$2$ is partially incorporated
in~\cite{Miranda.Morrison:1,Miranda.Morrison:2,Miranda.Morrison:book}:
for example, the isomorphism class of a finite quadratic form generated by an
element~$\Ga$ with $q(\Ga)=(u/p^k)\bmod\Z$,
which is $(2u/p^k)\bmod2\Z$ in our notation,
is designated by
the class of $2u$ in $(\Z_p\units)/(\Z_p\units)^2$.
\endwarning

Given
a lattice~$N$ and a prime $p\in\PP$, we define
the number $e_p:=e_p(N)\in\N$
and the subgroup
$\tSigma_p:=\tSigma_p(N)\subset\GG_0:=\{\pm1\}\times\{\pm1\}$
as in \eqref{eq.invariants}. Algorithms computing $e_p(N)$ and
$\tSigma_p(N)$ are given explicitly in ~\cite{Miranda.Morrison:2}.
Computations are in terms of $\rank N$, $\det N$, and $\dN$ only,
which means that the genus $g(N)$ determines $e_p(N)$,
$\tSigma_p(N)$ and, moreover, $\Coker\map$.
One has $e_p=1$ and $\tSigma_p=\GG_0$
for almost all $p\in\PP$.

\theorem[Miranda--Morrison~\cite{Miranda.Morrison:1,Miranda.Morrison:2}]\label{th.MM}
For $N$ as in~$\theN$,
there is an $\FF2$-module $\MM(N)$ and an exact sequence
\[*
\OG(N)\overset\map\longrightarrow\Aut\dN
\overset\mmh\longrightarrow\MM(N)\to\genus(N)\to\one,
\]
where $\genus(N)$ is the genus group of~$N$.
One has $\ls|\MM(N)|=e(N)/[\GG_0:\tSigma(N)]$, where
$e(N):=\prod_{p\in\PP}e_p(N)$ and
$\tSigma(N):=\bigcap_{p\in\PP}\tSigma_p(N)$.
\pni
\endtheorem

The group $\MM(N)$ and homomorphism
$\mmh\:\Aut\dN\to\MM(N)$ given by \autoref{th.MM}
will be called, respectively, the
\emph{Miranda--Morrison group} and \emph{Miranda--Morisson homomorphism}
of~$N$.
The next statement follows from \autoref{th.N.extensions}, \autoref{th.MM},
and the fact that a unimodular even indefinite lattice is unique in its
genus.

\corollary[Miranda--Morrison~\cite{Miranda.Morrison:1,Miranda.Morrison:2}]\label{cor.MM}
Let $L$ be a unimodular even lattice and $S\subset L$ a primitive sublattice
such that $N:=S^\perp$ is
as in~$\theN$.
Then the strict isomorphism classes of primitive extensions $S\pe L$ are
in a
canonical
one-to-one correspondence with the Miranda--Morrison group $\MM(N)$.
\pni
\endcorollary

Generalizing, fix an anti-isometry $\Gf\:\dS\to\dN$ and consider the induced
map $\map^\Gf\:\!\OG(S)\to\Aut N$, see \autoref{s.extensions}.
Since $\Im\map\subset\Aut\dN$ is a normal subgroup with abelian quotient,
this map factors to a
homomorphism
$\map^\perp\:\!\OG(S)\to\Aut\dN\to\MM(N)$ independent of~$\Gf$.
Then, the following statement is an immediate consequence of
Theorems~\ref{th.N.extensions} and~\ref{th.MM}.

\corollary\label{cor.O'}
Let $S\subset L$ be as in \autoref{cor.MM}, and let $\OG'\subset\OG(S)$ be a
subgroup. Then the $\OG'$-isomorphism classes of
primitive extensions $S\pe L$ are
in a
one-to-one correspondence with the
$\FF2$-module
$\MM(N)/\map^\perp(\OG')$.
\done
\endcorollary

\autoref{th.MM} and \autoref{cor.MM}
cover most of our needs. However, in a few special cases, we
need the more advanced treatment of~\cite{Miranda.Morrison:book}.
Introduce the groups
\[*
\GG_{p,0}:=\{\pm1\}\times\Z_p\units/(\Z_p\units)^2\subset
\GG_p:=\{\pm1\}\times\Q_p\units/(\Q_p\units)^2,\quad p\in\PP,
\]
and
\[*
\GAz:=\prod_p\GG_{p,0}\subset
\GA:=\GAz\cdot\sum_p\GG_p\subset
\GG:=\prod_p\GG_{p}.
\]
(Since the groups involved are multiplicative, although abelian, we
follow~\cite{Miranda.Morrison:book} and use~$\cdot$ to denote the sum of
subgroups. However, we retain the notation $\sum$ and $\prod$ to
distinguished between direct sums and products.
Thus, the adelic version $\GA$ is the set of sequences
$\{(s_p,t_p)\}\in\GG$ such that $(s_p,t_p)\in\GG_{p,0}$ for almost all~$p$.)
Let also $\GQ:=\{\pm1\}\times\Q\units\!/(\Q\units)^2\subset\GA$.
Then $\GAz\cdot\GQ=\GA$ and the intersection $\GAz\cap\GQ$ is the
group $\GG_0=\{\pm1\}\times\{\pm1\}$ introduced above.
We recall that
$\Q\units\!/(\Q\units)^2$ is the $\FF2$-module on the basis $\{-1\}\cup\PP$,
\ie, it is the set of all square free integers.

On various occasions we will also consider the following subgroups:
\roster*
\item
$\GG_p^{++}:=\{1\}\times\Z_p\units/(\Z_p\units)^2\subset\GG_{p,0}$;
\item
$\GG_{2,2}\subset\GG_2^{++}$ is the subgroup generated by $(1,5)$;
\item
$\GQ^{--}\subset\GQ$ is the subgroup generated by $(-1,-1)$ and $(1,p)$,
$p\in\PP$;
\item
$\GG_0^{--}:=\GQ^{--}\cap\GG_0\subset\GG_0$
is the subgroup generated by $(-1,-1)$.
\endroster

We denote by $\iota_p\:\GG_p\into\GA$, $p\in\PP$,
and $\iota_\Q\:\GQ\into\GA$ the inclusions.
The images
$\iota_\Q(1,q)$ and $\iota_q(1,q)$,
$q\in\PP$, differ by an element of $\prod_p\GG_p^{++}$, \viz., by the
sequence $\{(1,s_p)\}$, where $s_q=1$ and $s_p$ is the class of~$q$ in
$\Z_p\units/(\Z_p\units)^2$ for $p\ne q$.

Defined and computed in~\cite{Miranda.Morrison:book} are certain
$\FF2$-modules
\[*
\Sigma_p\dual(N):=\Sigma\dual(N\otimes\Z_p)\subset
\Sigma_p(N):=\Sigma(N\otimes\Z_p),
\]
which depend on the genus of~$N$ only.
One has $\Sigma_p\dual\subset\GG_{p,0}$, $\Sigma_p\subset\GG_p$, and
$\Sigma_p\subset\GG_{p,0}$ for almost all~$p$. (In fact, for almost all
$p\in\PP$ one has $\Sigma_p\dual=\Sigma_p=\GG_{p,0}$.)
Hence,
\[*
\Sigma\dual(N):=\prod_p\Sigma_p\dual(N)\subset\GAz,\qquad
\Sigma(N):=\prod_p\Sigma_p(N)\subset\GA.
\]
In these notations, the invariants used in \autoref{th.MM} are
\[
e_p(N)=[\GG_{p,0}:\Sigma_p\dual(N)],\qquad
\tSigma_p(N)=\Sigma_0\dual(N\otimes\Z_p):=\Gf_p\1(\Sigma_p\dual(N)),
\label{eq.invariants}
\]
where $\Gf_p\:\GG_0\to\GG_{p,0}$ is the projection, and $\MM(N)$ is the
quotient
$\GAz/\Sigma\dual(N)\cdot\GG_0$.
(Clearly, $\tSigma(N)=\Sigma\dual(N)\cap\GG_0$.)
Unfortunately, the map $\prod_p\Aut\dN_p\to \MM(N)$ given by \autoref{th.MM}
does not respect the product structures.
The following statement refines \autoref{th.MM}, separating the genus group
and the $p$-primary components.

\theorem[Miranda--Morrison~\cite{Miranda.Morrison:book}]\label{th.MM.refined}
Let~$N$ be
as in~$\theN$.
Then\rom:
\roster
\item\label{MM.genus}
there is an isomorphism
$\genus(N)=\GA/\Sigma(N)\cdot\GQ$
\rom(hence, $N$ is unique in its genus
if and only if $\GA=\Sigma(N)\cdot\GQ$\rom)\rom;
\item\label{MM.coker}
there is a commutative diagram
\[*
\CD
\Aut\dN=\prod_p\Aut\dN_p@>\mm>>\prod_p\Sigma_p(N)/\Sigma_p\dual(N)\\
@VVV@VV\Gb V\\
\Coker\map@>\cong>>\Sigma(N)/\Sigma\dual(N)\cdot(\Sigma(N)\cap\GQ),
\endCD
\]
where all maps are epimorphisms,
$\mm$ is the product of certain
epimorphisms $\mm_p\:\Aut\dN_p\onto\Sigma_p(N)/\Sigma_p\dual(N)$, $p\in\PP$,
and
$\Gb$
is the quotient projection.
\pni
\endroster
\endtheorem

\subsection{A few simple consequences}\label{s.MM.computation}
The homomorphism~$\mm$ in \autoref{th.MM.refined}\iref{MM.coker} is easily
computed on reflections: for a mirror $\xi\in\dN_r$, $r\in\PP$,
modulo $\Sigma_r\dual(N)$ one has
\[*
\mm_r(\refl\xi)=(-1,mr^k),
\quad\text{where}\quad
\xi^2=\frac{2m}{r^k}\bmod2\Z,\ \gcd(m,r)=1,\ k\in\N.
\]
If $r=2$ and $\xi^2=0\bmod\Z$, this value is only well defined
$\bmod\,\GG_2^{++}$; if $r=2$ and $\xi^2=\frac12\bmod\Z$, it is well defined
$\bmod\,\GG_{2,2}$. In these two cases, the disambiguation of
$\mm_r(\refl\xi)$
needs more information about~$\xi$ and~$N$:
one needs to represent~$\xi$ in the form
$\frac12x$ for some $x\in N\otimes\Z_2$.
Given another prime~$p$,
consider the homomorphism $\chi_p\:\Z_p\units/(\Z_p\units)^2\onto\{\pm1\}$,
\[*
\chi_p(m):=\Bigl(\frac{m}{p}\Bigr)\quad\text{if $p\ne2$},\qquad
\chi_2(m):=m\bmod4,
\]
and define the \emph{$p$-norm}~$\ls|\xi|_p\in\{\pm1\}$ and the
`Kronecker symbol'~$\Gd_p(\xi)\in\{\pm1\}$ \via
\[*
\ls|\xi|_p:=\begin{cases}
\chi_p(r^k),&\text{if $r\ne p$},\\
\chi_p(m),&\text{if $r=p$},
\end{cases}
\qquad
\Gd_p(\xi)=(-1)^{\Gd_{p,r}},
\]
where $\Gd_{p,r}$ is the conventional Kronecker symbol. (If $p=2$ and
$\xi^2=0\bmod\Z$, then $\ls|\xi|_2$ is undefined.)
Finally, introduce a few \latin{ad hoc}
notations for a lattice~$N$:
\roster*
\item
the group $\MM_p(N)=\{\pm1\}$ if $p=1\bmod4$ and
$e_p(N)\cdot\ls|\tSigma_p(N)|=8$; in all other cases,
$\MM_p(N)=\one$;
\item
the map~$\mma_p$ sending a mirror~$\xi$ to
$\ls|\xi|_p\in\MM_p(N)$, with the convention
that $\mma_p(\xi)=1$ whenever $\MM_p(N)=\one$;
\item
the map $\mmb_p$ sending a mirror~$\xi$ to an element of~$\GG_0$:
if $p=1\bmod4$, then
$\mmb_p(\xi)=(\Gd_p(\xi)\cdot\ls|\xi|_p,1)$;
otherwise, $\mmb_p(\xi)=\Gd_p(\xi)\times\ls|\xi|_p$.
\endroster

Following~\cite{Miranda.Morrison:book}, we say that a lattice~$N$ is
\emph{$p$-regular}, $p\in\PP$, if $\Sigma_p\dual(N)=\GG_{p,0}$, \ie, if
$e_p(N)=1$. We will also say that the prime~$p$ is \emph{regular} with respect to~$N$;
otherwise, $p$ is \emph{irregular}.
In several statements below,
we make a technical assumption that $\Sigma_2\dual(N)\supset\GG_{2,2}$;
this inclusion
does hold for the transcendental lattices of all primitive homological types
(see \autoref{s.htype})
except $\sset=\singset{A15+A3}$, see~\cite{Miranda.Morrison:book}.

\lemma\label{lem.p}
Let~$N$ be a lattice as in~$\theN$, $\Sigma_2\dual(N)\supset\GG_{2,2}$, and
assume that $N$ has one irregular prime~$p$.
Then $\MM(N)=\MM_p(N)$ and
$\mmh(\refl\xi)=\mma_p(\xi)$ for a mirror~$\xi$.
\endlemma

\lemma\label{lem.p,q}
Let~$N$ be a lattice as in~$\theN$, $\Sigma_2\dual(N)\supset\GG_{2,2}$, and
assume that $N$ has two irregular primes $p$, $q$.
Then
\[*
\MM(N)=\MM_p(N)\times\MM_q(N)\times(\GG_0/\tSigma_p(N)\cdot\tSigma_q(N))
\]
and one has
$\mmh(\refl\xi)=\mma_p(\xi)\times\mma_q(\xi)\times(\mmb_p(\xi)\cdot\mmb_q(\xi))$
for a mirror $\xi\in\dN$,
provided that $\xi^2\ne0\bmod\Z$ if $p=2$ or~$q=2$.
\endlemma

\corollary\label{cor.p}
Under the hypotheses of \autoref{lem.p,q}, assume, in addition, that
$\ls|\MM(N)|=\ls|\MM_p(N)|=2$.
Then $\MM(N)=\MM_p(N)$ and
$\mmh(\refl\xi)=\ls|\xi|_p$ for a mirror~$\xi$.
\done
\endcorollary

\proof[Proof of Lemmas~\ref{lem.p} and~\ref{lem.p,q}]
Let $\GG'_{p,0}:=\GG_{p,0}$ for $p\ne2$ and
$\GG'_{2,0}:=\GG_{2,0}/\GG_{2,2}$, so that
we can identify
$\GG'_{p,0}\cong\{\pm1\}\times\{\pm1\}$ for all $p\in\PP$.
If $p\ne1\bmod4$, the map $\Gf_p\:\GG_0\to\GG'_{p,0}$ is an epimorphism; if
$p=1\bmod4$, one has $\Gf_p(\GG_0)=\{\pm1\}\times\{1\}$.
Modulo~$\GQ^{--}$, the image $\mm(\refl\xi)$ equals
$\mma(\refl\xi):=\{(\Gd_s(\xi),\ls|\xi|_s)\}\in\prod\GG'_{s,0}$.

Now, the first statement of each lemma is a computation of the group
$\MM(N)=\GAz/\Sigma\dual(N)\cdot\GG_0$, which can be done in $\GG'_{p,0}$ or
$\GG'_{p,0}\times\GG'_{q,0}$;
our group $\MM_p(N)$ is the
quotient $\GG_{p,0}/\Sigma_p\dual(N)\cdot\Im\Gf_p$.
The second statement is the computation of the image of $\mma(\xi)$
in $\MM(N)$: the maps $\mma_p$ and $\mmb_p$ are the projections
$\GG_{p,0}\to\MM_p(N)$ and $\GG_{p,0}\to\Im\Gf_p$, respectively.
For the latter,
we use the following fact, see~\cite{Miranda.Morrison:book}: if
a prime $p=1\bmod4$
is irregular for~$N$ and $\Sigma_p\dual(N)\not\subset\Im\Gf_p$, then
$\Sigma_p\dual(N)$ is generated by $(-1,-1)$.
\endproof

\subsection{The positive sign structure}\label{s.signs}
A \emph{positive sign structure} on
a lattice~$N$ is a choice of an orientation of a maximal positive
definite subspace of $N\otimes\R$.
(Recall that the orthogonal projection of one such subspace to another is an
isomorphism and, hence, all these spaces admit a coherent orientation.)
We will use the map $\det_+\:\!\OG(N)\to\{\pm1\}$
sending an auto-isometry to $+1$ or $-1$ if it preserves or, respectively,
reverses a positive sign structure.
Thus,
$\Oplus(N):=\Ker\det_+$ is the subgroup
of auto-isometries preserving positive sign structures.
(In the notation
of~\cite{Miranda.Morrison:book}, one has $\det_+={\det}\cdot{\spin}$
and $\Oplus=\OG^{--}$.)
The following statement is essentially contained
in~\cite{Miranda.Morrison:book}.

\proposition[Miranda--Morrison~\cite{Miranda.Morrison:book}]\label{prop.MM.--}
Let~$N$ be a lattice as in~$\theN$.
Then
one has $\tSigma(N)\subset\GG_0^{--}$
if and only if $\det_+a=1$ for all
$a\in\Ker[\map\:\!\OG(N)\to\Aut\dN]$.
\pni
\endproposition

Thus, if $\tSigma(N)\subset\GG_0^{--}$, there is a well defined descent
$\det_+\:\Im\map\to\{\pm1\}$.
The next lemma computes the values of $\det_+$ on reflections.

\lemma\label{lem.det+}
Let~$N$ be a lattice as in~$\theN$, $\Sigma_2\dual(N)\supset\GG_{2,2}$, and
assume that there is a prime~$p$ such that
$\tSigma_p(N)\subset\GG_0^{--}$.
Then, for a mirror $\xi\in\dN$ such that $\refl\xi\in\Im\map$ and
$\xi^2\ne0\bmod\Z$ if $p=2$, one has
$\det_+\refl\xi=\Gd_p(\xi)\cdot\ls|\xi|_p$.
\endlemma

\proof
The proof is similar to that of Lemmas~\ref{lem.p} and~\ref{lem.p,q}:
we \emph{assume} that the element
$\mma(\refl\xi)\cdot\iota_\Q(\Gd_p(\xi),\Gd_p(\xi))$
representing~$\refl\xi$
lies in
$\Sigma\dual(N)\cdot\GG_0$ and compute its image in
$\Sigma\dual(N)\cdot\GG_0/\Sigma\dual(N)\cdot\GG_0^{--}=\{\pm1\}$.
This can be done in $\GG_{p,0}$.
\endproof

\autoref{prop.MM.--} can be restated in a form closer to
\autoref{th.MM}: introducing the group
$\MM_+(N):=\GAz/\Sigma\dual(N)\cdot\GG_0^{--}$, one has an exact sequence
\[
\OG_+(N)\overset\map\longrightarrow\Aut\dN
\overset{\mmh_+}\longrightarrow\MM_+(N)\to\genus(N)\to\one.
\label{eq.MM+}
\]
The groups $\MM_+(N)$, as well as a few other counterparts, are also computed
in~\cite{Miranda.Morrison:2}: for the order $\ls|\MM_+(N)|$, one merely
replaces $\tSigma(N)$ with $\tSigma(N)\cap\GG_0^{--}$ in \autoref{th.MM}.
In the special case of at most two irregular primes, the computation
is very similar to \autoref{s.MM.computation}.
For an irregular prime~$p$,
denote $\tSigma\plus_p(N):=\tSigma_p(N)\cap\GG_0^{--}\subset\GG_0^{--}$
and introduce
the groups $\MM\plus_p(N)$ and maps $\mma\plus_p$, $\mmb\plus_p$ defined on
the set of mirrors and
taking values in $\MM\plus_p(N)$ and $\GG_0^{--}=\{\pm1\}$, respectively,
as follows:
\roster*
\item
if $p=1\bmod4$, then $\MM\plus_p(N)=\MM_p(N)$, $\mma\plus_p=\mma_p$, and
$\mmb\plus_p(\xi)=\Gd_p(\xi)\cdot\ls|\xi|_p$;
\item
if $p\ne1\bmod4$, then $\MM\plus_p(N)=\GG_0/\tSigma_p(N)\cdot\GG_0^{--}$
(if $p\ne2$ or $\Sigma_2\dual(N)\supset\GG_{2,2}$, one has
$\MM\plus_p(N)=\{\pm1\}$ if $e_p(N)\cdot\ls|\tSigma\plus_p(N)|=4$
and $\MM\plus_p(N)=\one$ otherwise);
\item
if $p\ne1\bmod4$ and $\MM\plus_p(N)\ne\one$, then
$\mma\plus_p(\xi)=\Gd_p(\xi)\cdot\ls|\xi|_p$ and
$\mmb\plus_p(\xi)=\ls|\xi|_p$;
\item
if $p\ne1\bmod4$ and $\MM\plus_p(N)=\one$, then $\mma\plus_p(\xi)=1$ and
$\mmb\plus_p(\xi)$ is the image of $\mmb(\xi)=\Gd_p(\xi)\times\ls|\xi|_p$, see
\autoref{s.MM.computation}, under the projection
$\GG_0\to\GG_0/\tSigma_p(N)=\GG_0^{--}$.
\endroster
(In the last case, one has $\mmb\plus_p(\xi)=\ls|\xi|_p$ unless $p=2$.)
The proof of the next two statements repeats literally that of
Lemmas~\ref{lem.p} and~\ref{lem.p,q}.

\lemma\label{lem.p+}
Let~$N$ be a lattice as in~$\theN$, $\Sigma_2\dual(N)\supset\GG_{2,2}$, and
assume that $N$ has a single irregular prime~$p$.
Then one has $\MM_+(N)=\MM\plus_p(N)$ and
$\mmh_+(\refl\xi)=\mma\plus_p(\xi)$ for a mirror $\xi\in\dN$
such that $\xi^2\ne0\bmod\Z$ if $p=2$.
\done
\endlemma

\lemma\label{lem.p,q+}
Let~$N$ be a lattice as in~$\theN$, $\Sigma_2\dual(N)\supset\GG_{2,2}$, and
assume that $N$ has two irregular primes $p$, $q$.
Then
\[*
\MM_+(N)=\MM\plus_p(N)\times\MM\plus_q(N)\times
 (\GG_0^{--}\!/\tSigma\plus_p(N)\cdot\tSigma\plus_q(N))
\]
and
one has
$\mmh_+(\refl\xi)=\mma\plus_p(\xi)\times\mma\plus_q(\xi)\times(\mmb\plus_p(\xi)\cdot\mmb\plus_q(\xi))$
for a mirror $\xi\in\dN$ such that $\xi^2\ne0\bmod\Z$ if $p=2$ or~$q=2$.
\done
\endlemma

\corollary\label{cor.p+}
Under the hypotheses of \autoref{lem.p,q+}, assume, in addition, that
$\ls|\MM_+(N)|=\ls|\MM\plus_p(N)|=2$.
Then $\MM_+(N)=\MM\plus_p(N)$ and
$\mmh(\refl\xi)=\mma\plus_p(\xi)$ for a mirror
$\xi\in\dN$ such that $\xi^2\ne0\bmod\Z$ if $p=2$.
\done
\endcorollary

\section{The deformation classification}\label{S.classification}

\subsection{The homological type}\label{s.htype}
Consider a simple sextic $\sextic\subset\Cp2$.
Recall (see \autoref{s.sextics})
that we denote by $X\to\Cp2$ the minimal resolution of singularities
of the double covering of~$\Cp2$ ramified at~$\sextic$, and that the set of
singularities of~$\sextic$ can be identified with the sublattice
$\sset\subset\bL$ spanned by the classes of the exceptional divisors.
Let $\tau\: X\to X$
be the deck translation of the covering.

\lemma\label{lem.tau}
The induced action of~$\tau$ on the Dynkin graph $\Dynkin:=\Dynkin_\sset$
preserves the components of~$\Dynkin$\rom;
it acts by the only nontrivial symmetry on the components of type
$\bA_{p\ge2}$, $\bD\subtext{odd}$, or~$\bE_6$, and
by the identity otherwise.
\done
\endlemma

\remark\label{rem.tau}
In other words, $\tau\:\Dynkin\to\Dynkin$
can be characterized as the `simplest' symmetry of~$\Dynkin$ inducing
$-{\id}$ on $\discr\sset$.
\endremark

In addition to~$\sset$, we have the class $h\in\bL$ of the pull-back of a
generic line in~$\Cp2$. Obviously, $h$ is orthogonal to~$\sset$ and $h^2=2$.
Let $\sset_h:=\sset\oplus\Z h$. The triple
$\htype:=(\sset,h,\bL)$, \ie, the lattice extension $\sset_h\into\bL$
regarded up to isometries of~$\bL$ preserving~$\sset$ (as a set) and~$h$, is
called the \emph{homological type} of~$\sextic$. This extension is subject to
certain restrictions,
which are summarized in the following definitions.

\definition\label{def.htype}
Let $\sset$ be a root lattice.
A \emph{homological type} (extending~$\sset$)
is an extension $\sset_h:=\sset\oplus\Z h\into\bL$
satisfying the following conditions:
\roster
\item\label{htype.1}
any vector $v\in(\sset\otimes\Q)\cap\bL$ with $v^2=-2$ is in~$\sset$;
\item\label{htype.2}
there is no vector $v\in\tS:=(\sset_h\otimes\Q)\cap\bL$
with $v^2=0$ and $v\cdot h=1$.
\endroster
\enddefinition

Note that condition~\iref{htype.2} in this definition can be restated as
follows: if $a$ is a generator of an orthogonal summand $\bA_1\subset\sset$,
the vector $a+h$ is primitive in~$\bL$.

Given a homological type $\htype=(\sset,h,\bL)$, we let
\roster*
\item
$\tilde\sset:=(\sset\otimes\Q)\cap\bL$ be the primitive hull of~$\sset$,
\item
$\tS:=(\sset_h\otimes\Q)\cap\bL$ be the primitive hull of~$\sset_h$, and
\item
$\bT:=\sset_h^\perp$ with $\dT=\discr\bT$ be the \emph{transcendental
lattice}.
\endroster
Since $\sigma_+ \bT=2$, all positive definite $2$-spaces in $\bT\otimes\R$
can be oriented in a coherent way.
A choice $\pss$ of one of these coherent orientations,
\ie, a positive sign structure on~$\bT$, see \autoref{s.signs},
is called an \emph{orientation} of~$\htype$.
The homological type of a plane sextic~$\sextic$ has a canonical orientation,
\viz. the one given by the real and imaginary parts of the class of a
holomorphic form~$\Go$ on~$X$.


An \emph{automorphism} of a homological type $\htype=(\sset,h,\bL)$ is an
autoisometry
of $\bL$ preserving $\sset$ (as a set) and~$h$.
The group of automorphisms
of $\htype$ is denoted by $\Aut\htype$. Let $\Aut_+\htype\subset\Aut\htype$
be the subgroup of
the automorphisms inducing~$\id$ on~$\bT$. On the other hand,
we have the group $\Aut_h\tS\subset\OG(\tS)$ of the isometries of~$\tS$
preserving~$h$.
There are obvious homomorphisms
\[
\Aut_+\htype\into\Aut\htype\to\Aut_h\tS\into\OG(\sset),
\label{eq.Aut}
\]
where the latter inclusion is due
to \autoref{htype.1} in \autoref{def.htype},
as $\sset\subset\tS$ is recovered as
the sublattice generated by the roots orthogonal to~$h$.
If the primitive extension $\tS\pe\bL$
is defined by an anti-isometry $\Gf\:\discr\tS\to\dT$ (see
\autoref{s.extensions}),
so that we have
a homomorphism $\map^\Gf\:\Aut_h\tS\to\Aut\dT$,
then, for $\epsilon={+}$ or empty,
\[
\Im[\Aut_\epsilon\htype\to\Aut_h\tS]=(\map^\Gf)\1\map(\OG_\epsilon(\bT)).
\label{eq.Im.Aut}
\]

The deformation classification of sextics is based on the following
statement.

\theorem[see~\cite{degt:JAG}]\label{th.JAG}
The map sending a plane sextic $\sextic\subset\PP^2$ to
its oriented homological type establishes a bijection between the set of
equisingular deformation classes of simple sextics and the set of isomorphism
classes of oriented homological types.
Complex conjugate sextics have isomorphic homological types that differ by
the orientations.
\pni
\endtheorem

A homological type
is called \emph{symmetric} if it admits an orientation
reversing automorphism.
According to \autoref{th.JAG}, symmetric are the homological types
corresponding to \emph{real}, \ie, conjugation invariant components of
$\all(\sset)$.

Recall
that, in \autoref{s.sextics}, the equisingular strata $\all(\sset)$ were
subdivided into families $\all_*(\sset)$. The precise definition is as
follows: the subscript~$*$ is the sequence of invariant factors of the
kernel~$\CK$ of the finite index extension $\sset_h\fe\tS$.
(Obviously, $\CK$ is invariant under equisingular deformations.)
Theorems~\ref{th.special} and~\ref{th.1-torus} below single out the families
$\all_\ns$ and $\all_3$, which are of our primary interest; they correspond
to $\CK=0$ and $\CK=\CG3$,
respectively.


A homological type $\htype=(\sset,h,\bL)$ is called \emph{primitive} if
$\sset_h\subset\bL$ is a primitive sublattice, \ie, if $\CK=0$.
In this case, one has
$\discr\tS=\dS\oplus\<\frac12\>$ and
the
inclusion $\Aut_h\tS\into\OG(\sset)$, see~\eqref{eq.Aut}, is an isomorphism.

\theorem[see~\cite{degt:Oka}]\label{th.special}
A simple plane sextic~$\sextic$ is irreducible and non-special if and only if
its homological type is primitive.
\pni
\endtheorem

The fact that primitive homological types give rise
to irreducible sextics was also observed in \cite{Yang},
where the primitivity is stated as a sufficient condition.

\theorem[see~\cite{degt:Oka}]\label{th.1-torus}
A simple plane sextic~$\sextic$ is irreducible and \torus{p},
$p=1$, $4$, or~$12$,
if and only if the kernel~$\CK$ of the
extension $\sset_h\fe\tS$ is, respectively, $\CG3$,
$\CG3\oplus\CG3$, or $\CG3\oplus\CG3\oplus\CG3$.
\pni
\endtheorem

There is a similar characterization of other special sextics: a sextic is
irreducible and $\DG{2n}$-special, $n>3$,
if and only if the kernel~$\CK$
is
$\CG{n}$; one necessarily has $n=5$ or~$7$.
Note that these statements cover all
possibilities for the kernel~$\CK$ free of $2$-torsion, and $\CK$ has
$2$-torsion if and only if the sextic is reducible, see,
\eg,~\cite{degt:book}.

\subsection{Extending a fixed set of singularities $\sset$ to a sextic}\label{s.steps}
By
\autoref{th.JAG},
given a simple set of singularities~$\sset$, the connected components of the
space $\all(\sset)$ modulo
the complex conjugation $\conj\:\Cp2\to\Cp2$ are enumerated by the
isomorphism classes of the homological types extending~$\sset$.
If a subscript~$*$ is specified, the set $\pi_0(\all_*(\sset)/\!\conj)$
is enumerated by the extensions with the kernel~$\CK$
of the finite index extension $\sset_h\fe\tS$ in the given
isomorphism class.

We are interested in the sets of singularities~$\sset$ with
$\Gm(\sset)\le18$. In this case, $\bT$ is indefinite and $\rank\bT\ge3$;
hence, Miranda--Morrison's results apply and, with~$\CK$ and, hence, $\tS$
fixed, the further extensions $\tS\pe\bL$ are enumerated by the cokernel of
the well-defined homomorphism $\map^\perp\:\Aut_h\tS\to\MM(\bT)$, see
\autoref{s.MM}. In the special case $\CK=0$, due to the isomorphism
$\Aut_h\tS=\OG(\sset)$, we have a canonical bijection
\[
\pi_0(\all_\ns(\sset)/\!\conj)=\Coker[\map^\perp\:\OG(\sset)\to\MM(\bT)],
\label{eq.all_ns}
\]
assuming that $\sset_h$ does admit a primitive extension to~$\bL$
and taking for~$\bT$ any representative of the genus $\sset_h^\perp$.

\subsection{Proof of \autoref{th.classification}}\label{proof.classification}
By Theorems~\ref{th.JAG} and~\ref{th.special}, for the first part of the
statement it suffices to list (using \autoref{th.N.existence})
all sets of singularities extending to a primitive homological type; the
resulting list is compared against the list of all perturbations of the
maximizing sets obtained.
Since the homological type is primitive, we have
$\discr\tS=\dS\oplus\<\frac12\>$.

For the second part, let~$\sset$ be one of the sets of singularities found,
$\Gm(\sset)\le18$, and let $\bT$ be a representative of the genus
$\genus(\sset_h^\perp)$.
In most cases, \autoref{th.MM} gives us $\MM(\bT)=0$ and,
due to \autoref{cor.MM}, a primitive homological
type extending~$\sset$ is unique up to strict isomorphism.
In the remaining cases, it suffices to show that
the map $\map^\perp\:\!\OG(\bS)\to\MM(\bT)$
is onto, see~\eqref{eq.all_ns}.

\table
\caption{Exceptional sets of singularities (see \autoref{proof.classification})}\label{tab.Sigma}
\let\tabentryii\tabentry
\let\tabbox\tabboxii
\let\tabquad\quad
\tabs\tabii
E6+2A4+2A2//i.5\cr
A5+2A4+2A2+A1//i.5\cr
3A4+3A2//i.2x2\cr
\tabbreak
E7+A7+2A2//i.2,3\cr
E6+A7+A5//i.2,3\cr
2A7+2A2//i.2,3\cr
\tabbreak
A7+A5+A4+A2//i.2,3\cr
2A6+2A2+2A1//i.3,7\cr
2A9//i.2A9\cr
\endtab\endtabs

\endtable

There are $32$ sets of singularities containing a point of type~$\bA_4$ and
satisfying the hypotheses of \autoref{lem.p} or \autoref{cor.p}
(with $p=5$); in these cases, a nontrivial symmetry of any type~$\bA_4$
points maps to the generator $-1\in\MM(\bT)$.
The remaining nine sets of singularities are
collected in \autoref{tab.Sigma}, with references to the list below, where we
indicate the Miranda--Morrison homomorphism
$\mmh\:\Aut\dT\to\MM(\bT)$ (given by \autoref{lem.p,q}) and automorphism(s)
of~$\sset$ generating $\MM(\bT)$.
\roster
\item\label{i.5}
$\mmh\:\refl\xi\mapsto\Gd_3(\xi)\cdot\Gd_5(\xi)\cdot\ls|\xi|_5\in\{\pm1\}$;
a transposition $\bA_4\leftrightarrow\bA_4$;
\item\label{i.2x2}
$\mmh\:\refl\xi\mapsto(\Gd_3(\xi)\cdot\Gd_5(\xi)\cdot\ls|\xi|_5,\ls|\xi|_5)\in\{\pm1\}\times\{\pm1\}$;
a symmetry of~$\bA_4$
\emph{and}
a transposition $\bA_4\leftrightarrow\bA_4$ (two generators);
\item\label{i.2,3}
$\mmh\:\refl\xi\mapsto\Gd_2(\xi)\cdot\Gd_3(\xi)\cdot\ls|\xi|_2\cdot\ls|\xi|_3\in\{\pm1\}$;
a transposition $\bA_2\leftrightarrow\bA_2$
or
a symmetry of~$\bA_4$, $\bA_5$, or~$\bE_6$;
\item\label{i.3,7}
$\mmh\:\refl\xi\mapsto\Gd_3(\xi)\cdot\Gd_7(\xi)\cdot\ls|\xi|_3\cdot\ls|\xi|_7\in\{\pm1\}$;
a transposition $\bA_1\leftrightarrow\bA_1$;
\item\label{i.2A9}
$\mmh\:\refl\xi\mapsto\ls|\xi|_5\in\{\pm1\}$; none.
\endroster
The last case $\sset=2\bA_9$ is special: the map
$\map^\perp\:\!\OG(\sset)\to\MM(\bT)$ is not surjective and there are two
deformation families, as stated.

To complete the proof,
we need to analyze whether the space
$\all_\ns(\sset)$ contains a real curve and, if it does not, whether
the homological type $\htype$ extending~$\sset$
is symmetric. This is done in \autoref{s.real.ns} below.
\qed

\subsection{Proof of \autoref{cor.degeneration}}\label{proof.degeneration}
Unless $\sset=\singset{2A9}$, the statement follows immediately from
\autoref{th.classification}. Indeed, there is a degeneration $\sset\dg\sset'$
to a maximizing set of singularities~$\sset'$.
Due to \cite[Proposition 5.1.1]{degt:8a2}, there is a degeneration
$\sextic\dg\sextic'$ of \emph{some} sextics $\sextic\in\all_\ns(\sset)$ and
$\sextic'\in\all_\ns(\sset')$. Since $\all_\ns(\sset)/\conj$ is connected, a
degeneration exists for \emph{any} sextic $\sextic\in\all_\ns(\sset)$.
The exceptional case $\sset=\singset{2A9}$
with disconnected moduli space
is given by \autoref{prop.2A9}, see
\autoref{proof.2A9} below.
\qed

\subsection{Proof of \autoref{prop.2A9}}\label{proof.2A9}
For $\sset_0=\singset{2A9}$, one has
$\bT\cong\Z u\oplus\Z v\oplus\Z w$, with $u^2=v^2=10$, $w^2=-2$.
The group~$\dT$ is
$\<\frac25\>\oplus\<\frac25\>\oplus\<\frac12\>\oplus\<\frac12\>\oplus\<\frac32\>$,
and $\Aut\dT$ is generated by
\[*
\Gs_{1,2}\:\Ga_{1,2}\mapsto-\Ga_{1,2},\quad
\Gs_3\:\Ga_1\leftrightarrow\Ga_2,\quad
\Gs_4\:\Ga_3\leftrightarrow\Ga_4.
\]
Let $\dS_h:=\discr\tS=\dS_0\oplus\<\frac12\>$.
According to \autoref{s.root},
the image of $\map\:\OG(\sset_0)\to\Aut\dS_h$ is generated by $-{\id}$ on each
of the two copies of $\discr\singset{A9}$
and by the transposition of the two copies.
Since $\ls|\MM(\bT)|=2$, the image
$\Im[\map\:\OG(\bT)\to\Aut\dT]$
is generated by the images $\Gs_1$, $\Gs_2$,
$\Gs_3\Gs_4$
of the auto-isometries
$u\mapsto -u$, $v\mapsto -v$, $u\leftrightarrow v$, respectively.
It is straightforward that
$\Im\map^\perp=0\subset\MM(\bT)$;
hence, by \autoref{cor.O'}, $\singset{2A9}\oplus\Z h$
extends to $\bL$ in two ways.
The proof of the fact that
both homological types are represented by real curves
is postponed till
\autoref{s.real}
below.

The two homological types can be distinguished as follows.
In
$\dT$,
there are
two non-characteristic elements of square~$\frac12$
and two pairs of opposite elements of
square $\frac25$,
and
the map
$\frac12u\mapsto\frac15u$,
$\frac12v\mapsto\pm\frac15v$
establishes a
bijection between these two-element sets.
A similar bijection in the other group $\dS_h$ is due to the
decomposition
$\dS_h=2\discr\bA_9\oplus\<\frac12\>$.
The two homological types extending~$\singset{2A9}$
differ by whether the anti-isometry
$\dS_h\to\dT$ does or does not respect these bijections.

Now, a simple computation shows
that each of the two sublattices $\sset_0\oplus\Z h\subset\bL$ extends to
both $\sset_i\oplus\Z h\subset\bL$, $i=1,2$
(where $\sset_1=\singset{A19}$ and
$\sset_2=\singset{A10+A9}$ are as in the statement),
and these are all possible
degenerations of~$\sset_0$.
On the other hand, each~$\sset_i$, $i=1,2$,
extends to two distinct real homological types,
see~\cite{Shimada:maximal}, and each
of the resulting families admits a unique, up to deformation,
perturbation to \singset{2A9}, \cf. \cite[Proposition 5.1.1]{degt:8a2}.
These observations complete the proof.
\qed

\subsection{Proof of \autoref{th.torus} and \autoref{cor.torus}}\label{proof.torus}
Let $\sset$ be a set of singularities of weight~$6$ or~$7$. As shown
in~\cite{degt:Oka}, up to automorphism of~$\sset$, there is at most one
isotropic order~$3$ element $\Gb\in\dS$ satisfying condition~\iref{htype.1}
in \autoref{def.htype}. Such an element does exist if and only if
$w(\sset)=6$ or $w(\sset)=7$ and $\sset$ contains $\bA_2$ as a direct
summand. (In the latter case, the extra $\bA_2$ point becomes an outer
singularity; all other singular points of positive weight are inner.)
This element $\Gb$ has the form $\sum_i(\pm\Ga_i)$, where
$\Ga_i$ are the only (up to sign) order~$3$ elements in the discriminants of
the inner singular points. Important for
Theorems~\ref{th.N.existence} and~\ref{th.MM}
is the relation
between~$\dS$ and $\tilde\dS:=\discr\tilde\sset$. One has:
\roster*
\item
$\ell_p(\tilde\dS)=\ell_p(\dS)$ and $\det_p\tilde\dS=\det_p\dS$ for all
primes $p\ne3$;
\item
$\ls|\tilde\dS|=\frac19\ls|\dS|$ and $\det_3\tilde\dS=-9\det_3\dS$;
\item
$\ell_3(\tilde\dS)=\ell_3(\dS)-\Gd$, where
$\Gd=1$ if $\sset$ contains (as a direct summand)
$\bA_{17}$ or $2\bA_8$ and $\Gd=2$ otherwise.
\endroster

Now, as in \autoref{proof.classification}, we compare two lists:
the sets of singularities extending to a homological types with kernel $\CG3$
(using \autoref{th.N.existence}) and those obtained by perturbations from the
maximizing sets, see \autoref{tab.torus}. These lists coincide.
For each set of singularities~$\sset$ found, \autoref{th.MM} gives us
$\MM(\bT)=0$; hence, there is a unique homological type and
the space $\all_3(\sset)/\conj$ is connected.
In view of the first part, this fact implies \autoref{cor.torus},
and it remains
to analyze the real structures. This is done in \autoref{s.real.torus} below.
\qed

\subsection{Digression: permutations of the singular points}\label{s.permutations}
Consider a sextic~$\sextic$ with the set of singularities~$\sset$, and let
$\all(\sextic)$ be the \emph{connected} equisingular stratum
containing~$\sextic$.
Denoting by $\SSG(\sset)$ the group of the
type-preserving
permutations of the
singular points constituting~$\sset$,
we have the so-called
\emph{monodromy representation} $\pi_1(\all(\sextic))\to\SSG(\sset)$.
In this section, we are interested in the image $\SSGp(\sextic)$ of this
homomorphism.
In other words, we can consider the covering
$\tilde\all(\sextic)\to\all(\sextic)$ whose points are sextics with marked
singular points; then, $[\SSG(\sset):\SSGp(\sextic)]$ is the
number of the connected components of $\tilde\all(\sextic)$.

\theorem\label{th.group}
The permutation group $\SSGp:=\SSGp(\sextic)$ of a
non-special irreducible simple sextic~$\sextic$
with the set of singularities~$\sset$ is as follows\rom:
\roster*
\item
if $\Gm(\sextic)=19$, then $\SSGp$ is the group of permutations of the
$\bE_8$ points of~$\sset$\rom;
\item
if $\sset$ is one of the sets of singularities listed in \autoref{tab.group},
then $\SSGp$ is as shown in the table \rom(see the explanation after the
statement\rom).
\endroster
In all other cases, one has $\SSGp=\SSG(\sset)$.
\endtheorem

\table
\caption{Permutation groups (see \autoref{th.group})}\label{tab.group}
\tabs\tabii
[3E6]\cr
[2E6]+A6\cr
[2E6]+A5+A1\cr
[2E6]+A5\cr
[2A2]+E7+A7\cr
\tabbreak
[2A7+2A2]\cr
[2A2+2A1]+2A6\cr
[2A4]+E6+2A2\cr
[2A4]+A5+2A2+A1\cr
\tabbreak
[2A4+2A2]+D6\cr
[2A4+2A2]+2A3\cr
[2A4+3A2]+A3+A1\cr
[3A4]+[3A2]\cr
\endtab\endtabs

\endtable

The groups $\SSGp(\sextic)$ are encoded in \autoref{tab.group} by means of
one or several subsets $\sset_1,\sset_2,\ldots$ enclosed in brackets: a
permutation $\Gs\in\SSG(\sset)$ belongs to $\SSGp(\sextic)$ if and only if
the restriction of~$\Gs$ to each subset~$\sset_i$ is even.
Note that, in many cases, this condition actually implies that
$\SSGp(\sextic)$ is the trivial group.

\proof
If $\Gm(\sextic)=19$, then $\SSGp(\sextic)$ is the group of projective
symmetries of~$\sextic$; these groups are described in~\cite{degt:symmetric}.

In general, let $(\htype,\pss)$ be the oriented homological type
of~$\sextic$. From the description of the equisingular moduli spaces of
sextics, see, \eg,~\cite{degt:JAG}, it is immediate that the monodromy
representation can be factored as
\[*
\pi_1(\all(\sextic))\onto\Aut_+\htype\to\OG(\sset)\onto\SSG(\sset),
\]
where the arrow in the middle is the homomorphism~\eqref{eq.Aut}.
If $\htype$ is primitive and $\Gm(\sset)\le18$, we have a
well-defined homomorphism $\map^\perp\:\!\OG(\sset)\to\MM_+(\bT)$,
\cf. \autoref{s.MM}, where $\bT$ is the transcendental lattice;
this homomorphism factors through $\map'\:\Sym'\Dynkin_\sset\to\MM_+(\bT)$,
see \autoref{lem.root}.
Hence, combining the above observation with~\eqref{eq.MM+}
and~\eqref{eq.Im.Aut}, we conclude that
$\SSGp\subset\SSG(\sset)$ is the image of $\Ker\map'$.

The groups $\MM_+(\bT)$ are computed using
Lemmas~\ref{lem.p+} and~\ref{lem.p,q+}. For most curves, one has
$\MM_+(\bT)=\one$ and hence $\SSGp=\SSG(\sset)$.

There are $171$ sets of singularities~$\sset$ containing a point of
type~$\bA_2$ and satisfying the hypotheses of \autoref{lem.p+} or
\autoref{cor.p+} with $p=3$. For such curves, a non-trivial symmetry
of~$\bA_2$ maps to the generator $-1\in\MM_+(\bT)$; hence,
$\SSGp=\SSG(\sset)$.

Similarly, there are $28$ sets of singularities~$\sset$ containing a point of
type~$\bA_4$ and satisfying the hypotheses of \autoref{lem.p+} or
\autoref{cor.p+} with $p=5$: a non-trivial symmetry
of~$\bA_4$ maps to the generator $-1\in\MM_+(\bT)$.

In the very few remaining cases, the group $\Sym'\Dynkin_\sset$, identified with
its image in $\Aut\dS$, see \autoref{lem.root}, is generated by reflections,
and the map $\map'$ is computed explicitly using Lemmas~\ref{lem.p+}
and~\ref{lem.p,q+}. Details are left to the reader.
\endproof

\section{The fundamental group}\label{S.pi1}

\subsection{The degeneration principle}
Our computation of the fundamental groups is indirect;
it is based on a few previously known results and the following statement,
often referred to as the \emph{degeneration principle}.

\theorem[Zariski~\cite{Zariski:group}]\label{th.Zariski}
If a plane curve~$\sextic'$ degenerates to a reduced plane curve~$\sextic$,
there is an epimorphism
$\pi_1(\Cp2\sminus\sextic)\onto\pi_1(\Cp2\sminus\sextic')$.
\pni
\endtheorem

\corollary\label{cor.abelian}
If a plane sextic~$\sextic'$ degenerates to~$\sextic$ and
$\pi_1(\Cp2\sminus\sextic)=\CG6$, then also $\pi_1(\Cp2\sminus\sextic')=\CG6$.
\done
\endcorollary

\corollary\label{cor.MG}
If a sextic~$\sextic'$ of torus type degenerates to~$\sextic$ and
$\pi_1(\Cp2\sminus\sextic)=\MG$, then also $\pi_1(\Cp2\sminus\sextic')=\MG$.
\endcorollary

\proof
Since any sextic~$\sextic'$ of torus type is a degeneration of Zariski's
six-cuspidal sextic, there is an epimorphism
$\pi_1(\Cp2\sminus\sextic')\onto\MG$,
see~\cite{Zariski:group} and \autoref{th.Zariski}.
Since $\MG$ is a Hopfian group,
the statement follows from \autoref{th.Zariski}.
\endproof

\subsection{Proof of \autoref{cor.pi1}}\label{proof.pi1}
We need a slightly stronger statement, which is proved in the same way as
\autoref{cor.degeneration}, see \autoref{proof.degeneration}, by comparing
two independent lists: with few exceptions listed below, any non-special
irreducible plane sextic degenerates to one with known \emph{abelian} fundamental
group.

The exceptions are the six sets of singularities listed in
\autoref{th.classification} and
\[*
\gathered
\singset{2A4+2A3+2A2}\dg\singset{E8+A4+A3+2A2},\\
\singset{3A4+3A2},\ \singset{2A4+A3+3A2+A1}\dg\singset{E7+2A4+2A2}.
\endgathered
\]
The fundamental groups of the curves listed in \autoref{th.classification}
are computed in~\cite{degt:book}, using the degenerations
\[*
\singset{2D9}\dg\singset{D10+D9},\quad
\singset{2D7+2A2}\dg\singset{D10+D7+A2}
\]
to reducible maximizing sextics.
The groups of \emph{some} curves realizing the three other sets of
singularities are computed together with those of the corresponding
maximizing sextics, by analyzing the perturbations (see~\cite{degt:book} for
references). In view of the uniqueness given by \autoref{th.classification},
the results hold for \emph{all} curves.
\qed

\subsection{Proof of \autoref{cor.pi1.torus}}\label{proof.pi1.torus}
With one exception, \viz. the set of singularities \singset{(A8+A5+A2)+A4}, the
fundamental groups of all maximizing irreducible sextics of torus type are
known, see \cite{degt:book,degt:tetra} for references.
Comparing the two lists, one can easily see that all
but $14$ non-maximizing deformation families degenerate to
maximizing sextics~$\sextic$
with $\pi_1(\Cp2\sminus\sextic)=\MG$ known; for these curves, the fundamental
group is~$\MG$ due to \autoref{cor.MG}.
All sextics with at least one type $\bE_6$ type point are treated
in~\cite{degt:book}. The remaining exceptions are
\[*
\singset{(6A2)+4A1}\dg\singset{(6A2)+A3+2A1},
\]
studied in~\cite{degt:e6} as perturbations of \singset{(3E6)+A1}, and
\[*
\singset{(6A2)+A4+A1}\dg\singset{(A5+4A2)+A4+A1},
\]
studied in~\cite{degt:tetra} as perturbations of
\singset{(A8+3A2)+A4+A1}.
\qed

\section{Real structures}\label{S.real}

\subsection{Real sextics}\label{s.real}
A \emph{real structure} on a complex analytic variety $X$ is an anti-holomorphic
involution  $c\:X\rightarrow X$.
A \emph{real variety} is a pair $(X,c)$, where $X$ is a complex variety
and
$c$
is a real structure.
The fixed point set $X_{\R}:=\Fix c$
is called the \emph{real part}
of~$X$. (We routinely omit~$c$ in the notation when it is understood.)

Let $(X,c)$ be a real surface. A curve $\sextic\subset X$
is said to be \emph{real} if $c(D)=D$.
If $\bar{X}\to X$
is a double covering
branched over a (nonempty) real curve,
the real structure~$c$ lifts to two distinct real structures on $\bar{X}$;
the two lifts differ by the deck translation of the covering,
and all three involutions commute.

Any real structure on $\Cp2$ is equivalent to the standard complex
conjugation; in appropriate homogeneous coordinates, it is given by
$(z_0:z_1:z_2)\mapsto(\bar z_0:\bar z_1:\bar z_2)$.
In these coordinates, real curves are those defined by real polynomials.

\theorem\label{th.real}
A homological type~$\htype$ is realized by a real sextic if and only if
$\htype$ admits an \emph{involutive}
orientation reversing automorphism.
\endtheorem

\proof
The necessity is obvious: the real structure on~$\Cp2$ lifts to a real
structure on the covering $K3$-surface~$X$, which induces an involutive
automorphism of the homological type.

For the converse,
let $a\in\Aut\htype$ be an automorphism as in the statement.
Due to \autoref{lem.root}, the restriction $a|_\sset$ has the form
$r\circ(-s_*)$, where $r\in\Ker\map$ and $s_*$ is induced by an
\emph{involutive} symmetry $s\in\Sym'\Dynkin_\sset$.
Since
$\Ker\map\subset\Aut\htype$ (in the obvious way: automorphisms extend to
$\sset^\perp$ by the identity, see \autoref{th.N.isometry}),
the involution $r\1\circ a$ is also in
$\Aut\htype$. Let $c:=r\1\circ a\circ\refl{h}\in\OG(\bL)$; it is still an
involution and $c|_\bT=a|_\bT$.

Let $T_\pm$ be the $(\pm1)$-eigenspaces of the action of~$c$
on~$\bT\otimes\R$. Since $c$ reverses the orientation, one has
$\Gs_+T_\pm=1$. Hence, one can choose
generic
(\ie, maximally irrational)
vectors $\Go_\pm\in T_\pm$
such that $\Go_+^2=\Go_-^2>0$ and take $\Go:=\Go_++i\Go_-$ for the
class of a holomorphic form.
Let, further, $S_-$ be the $(-1)$-eigenspace of the action of~$c$ on
$\tS\otimes\R$. Since $h\in S_-$, one has $\Gs_+S_-=1$.
By the construction, $-c$ preserves a Weyl chamber of~$\sset$; hence,
condition~\iref{htype.1} in \autoref{def.htype} implies that $S_-$ is
\emph{not} orthogonal to a vector $v\in\tS$ of square~$(-2)$
and one can find a generic vector $\Gr\in S_-$, $\Gr^2>0$, and take it for
the class of a K\"{a}hler form.
These choices define a $2$-polarized $K3$-surface~$X$
with $\QOPNAME{Pic}X=\tS$
and, by an equivariant version of the
global Torelli theorem, $c$ is induced by a
real structure on~$X$ commuting with
the deck translation~$\tau$ of the ramified covering $X\to\Cp2$ defined
by~$h$.
This real structure descends to~$\Cp2$ and makes the sextic corresponding
to~$X$ (\ie, the branch curve) real.
\endproof

Let $\sextic$ be a real sextic with the set of singularities~$\sset$. The
real structure~$c$ lifts to two real structures on the covering $K3$-surface;
they take exceptional divisors to exceptional divisors and, hence, induce two
involutive
symmetries $c_\pm\:\Dynkin\to\Dynkin$ of the Dynkin graph
$\Dynkin:=\Dynkin_\sset$.
Define another symmetry $c_0\:\Dynkin\to\Dynkin$ as follows: on each
connected component $\Dynkin_i$ of~$\Dynkin$ \emph{fixed} by~$c_\pm$ and of
type other than~$\bD\subtext{even}$ let $c_0=\id$; on all other components,
let $c_0=c_\pm$.
In other words, since $c_-=c_+\circ\tau$, we just let $v\ra c_0=v$ for each
vertex~$v$ such that $v\ra c_+\ne v\ra c_-$, see \autoref{lem.tau}.

\corollary\label{cor.real.pert}
If a homological type~$\htype$ is realized by a real sextic $(\sextic,c)$,
then any $c_0$-invariant perturbation
$\htype'$ of~$\htype$ is also realized by a real
sextic~$\sextic'$.
\endcorollary

Note that we do \emph{not} assert that $\sextic'$ degenerates
to~$\sextic$ in the class of real sextics. A real perturbation can be
found if $\htype'$ is invariant under one of $c_\pm$.

\proof[Proof of \autoref{cor.real.pert}]
Let $c_*\:\bL\to\bL$ be the automorphism of~$\htype$
induced by one of the two lifts
of~$c$. Composing~$c_*$ with~$-\tau_*$ on \emph{some} of the indecomposable
summands of~$\sset$, we can change it to another involutive
automorphism~$c'$ of $\htype$ (see \autoref{lem.root} and
\autoref{th.N.isometry})
inducing~$c_0$ on~$\Dynkin$.
Then $c'$ preserves~$\sset'$; hence, $c'\circ\refl{h}$
can be regarded as an involutive
orientation reversing automorphism of~$\htype'$, and \autoref{th.real}
applies.
\endproof

\subsection{End of the proof of \autoref{th.classification}}\label{s.real.ns}
It is easily confirmed that
most sets of singularities~$\sset$ with $\Gm(\sset)\le18$
are \emph{symmetric} perturbations of maximizing
sets of singularities realized by real sextics, see Tables~\ref{tab.triple}
and~\ref{tab.double}. (In the tables, marked with a $^*$ are pairs of
isomorphic singular points permuted by the complex conjugation. These pairs
should be taken into account when analyzing symmetric perturbations. Note
that singular points of type $\bD\subtext{even}$ do not appear in irreducible
maximizing sextics.) Due to \autoref{cor.real.pert}, these sets of
singularities are realized by real curves.

The remaining $25$ sets of singularities
are listed in \autoref{tab.nonreal}. Each of these sets~$\sset$ extends to a
unique
(up to isomorphism) primitive homological type~$\htype$,
and we denote by~$\bT$ the
corresponding transcendental lattice.
In each case, the natural
homomorphism $\map\:\!\OG(\bT)\to\Aut\dT$ is surjective.

\table
\caption{Exceptional sets of singularities}\label{tab.nonreal}
\def\*{\llap{$\specialmark$\,}}
\def\?{\llap{$^?$}}
\def\tabentryii#1\end{\tabboxii{\quad\tabentryiii#1//}}
\def\tabentryiii#1//#2//{#2\singset{#1}}
\tabs\tabii
[3A6]_7//\cr
[2A6]_7+D6//\cr
[2A6]_7+D5+A1//\cr
[2A6]_7+2A3//\cr
[2A5]_3+E8//\cr
[E6+A11]_3+A1//\cr
[E6+A5]_3+E7//\cr
[E6+A5]_3+A7//\cr
\tabbreak
[E6+A5]_3+A6+A1//\cr
[E6+2A5]_3+A1//\cr
[E7+A7]_2+A4//\cr
[2A7]_2+A4//\cr
A7+A6+A5//\cr
2D7+2A2//\cr
D7+D4+3A2//\cr
2D4+4A2//\cr
\tabbreak
2E7+A4//\*\cr
E7+D5+A6//\*\cr
E7+A11//\*\cr
E7+A6+A5//\*\cr
2D9//\*\cr
D9+D8//\*\cr
2D8//\*\cr
D5+A7+A6//\*\cr
E7+2A4+A3//\cr
\endtab\endtabs

\endtable

By \autoref{th.N.isometry}, the homological type~$\htype$ is symmetric if and
only if there is an isometry $a\in\OG(\bT)$ with $\det_+a=-1$ and such that
$\map(a)\in\map^\Gf(\OG(\sset))$, where $\map^\Gf$ is induced by any
anti-isometry $\Gf\:\dS\oplus\<\frac12\>\to\dT$.
If (and only if) $a$ as above can be chosen \emph{involutive}, then so is
$\map(a)$ and, due to \autoref{lem.root}, $a$ extends to~$\bL$
by an \emph{involutive} isometry of~$\sset$;
hence,
$\all_\ns(\sset)$ contains real curves, see \autoref{th.real}.

\lemma\label{lem.asymmetric}
The first twelve sets of singularities in \autoref{tab.nonreal} \rom(those with
a $[\,\cdot\,]_p$ pattern\rom) extend to \emph{asymmetric} primitive
homological types.
\endlemma

\proof
Let $\sset$ be one of the sets of singularities in question. Then
$\tSigma(\bT)\subset\GG_0^{--}$, see \autoref{s.signs}, and there is a
well defined map $\det_+\:\Aut\dT\to\{\pm1\}$. We use \autoref{lem.det+}
(with the `test prime'~$p$ indicated in the table) to
show that $\det_+$ takes value~$+1$ on the image of $\OG(\sset)$.
If $p=7$ (the first four lines), the latter image is generated by reflections
$\refl\xi$ such that either
\roster*
\item
$\xi^2=\frac67$ (a symmetry of the Dynkin graph of~$\bA_6$), or
\item
$\xi^2=\frac{12}7$ (interchanging of two copies of~$\bA_6$), or
\item
$\xi\in\dT_2$ (isometries involving the other singular points);
\endroster
on the other hand, one has $(\frac{-3}7)=(\frac{-6}7)=(\frac27)=1$.
If $p=3$ (the next six sets of singularities), the image of $\OG(\sset)$ is
generated by the following automorphisms $a$:
\roster*
\item
$\refl\xi$ with $\xi^2=\frac43$
(a symmetry of the Dynkin graph of~$\bE_6$ or~$\bA_5$),
\item
$\refl\xi\refl\eta$ with $\xi^2=\frac23$, $\eta^2=1$
(interchanging of two copies of $\bA_5$),
\item
$\refl\xi\refl\eta$ with $\xi^2=\frac23$, $\eta^2=\frac14$
(a symmetry of the Dynkin graph of~$\bA_{11}$),
\item
$\refl\xi$ with $\xi^2=\frac78$ or $\xi^2=\frac67$
(a symmetry of the Dynkin graph of~$\bA_7$ or~$\bA_6$).
\endroster
In each case, \autoref{lem.det+} (with $p=3$) implies that $\det_+a=1$.
Finally, if $p=2$ (the last two sets of singularities),
we have reflections $\refl\xi$ such that either
\roster*
\item
$\xi^2=\frac78$ (a symmetry of the Dynkin graph of~$\bA_7$), or
\item
$\xi^2=\frac74$ (interchanging of two copies of~$\bA_7$), or
\item
$\xi^2=\frac45$ (a symmetry of the Dynkin graph of~$\bA_4$).
\endroster
\autoref{lem.det+} (with $p=2$)
implies that $\det_+\refl\xi=1$.
\endproof

Listed in the last column in \autoref{tab.nonreal} are the sets of
singularities~$\sset$
extending to symmetric homological types due to \autoref{prop.MM.--}.
However, since we want to represent these types by real
sextics, we will attempt to find \emph{involutive} orientation reversing
automorphisms, see \autoref{th.real}. A simplest automorphism with this
property would be a reflection $\refl{a}$, $a\in\bT$, $a^2=2$.

\lemma
If $\sset$ is one of the sets of singularities marked with a~$\specialmark$ in
\autoref{tab.nonreal}, the lattice~$\bT$ contains a vector~$a$ with $a^2=2$.
\endlemma

\proof
It suffices to find an embedding $\sset_h\oplus\Z a\into\bL$,
$a^2=2$, with the image of~$\sset_h$ primitive.
In each case, there is an element $\Ga\in\discr\sset_h$ with
$\Ga^2=-\frac12\bmod2\Z$. Let $\Gb\in\discr(\Z a)=\<\frac12\>$ be the generator,
and let
$\sset_h'$ be the finite index extension of~$\sset_h$
with the kernel generated by
$\Ga+\Gb$.
On a case-by-case basis one confirms that \autoref{th.N.existence} implies
the existence of a primitive embedding $\sset_h'\into\bL$.
(In the last case, the set of singularities \singset{D5+A7+A6}, the
element~$\Ga$ above should be chosen carefully,
\viz. $\Ga=2\Ga_1+4\Ga_2+\Ga_4$ in
$\discr\sset_h=\<\frac34\>\oplus\<\frac98\>\oplus\<\frac87\>\oplus\<\frac12\>$.)
\endproof

The set of singularities \singset{A7+A6+A5} is considered in
\autoref{prop.A7+A6+A5}, see \autoref{proof.A7+A6+A5} below, and
the remaining four deformation families are real and contain real curves; for
proof, we construct explicit reflections in $\OG(\bT)$.

If $\sset=\singset{2D7+2A2}$, then $\bT=\Z u\oplus\Z v\oplus\Z w$ with
$u^2=4$, $v^2=-12$, $w^2=6$,
and the reflection $\refl{u}$ extends to an involutive automorphism of~$\htype$
(\via\ $-{\id}$ on one of the $\bD_7$ components).
Hence, $\all_\ns(\sset)$ contains a real curve; by \autoref{cor.real.pert},
so do $\all_\ns(\singset{D7+D4+3A2})$ and $\all_\ns(\singset{2D4+4A2})$.

Finally,
if $\sset=\singset{E7+2A4+A3}$, then $\bT=\Z u\oplus\Z v\oplus\Z w$ with
$u^2=v^2=10$, $w^2=-4$.
Since $\map\:\!\OG(2\bA_4)\to\discr2\bA_4$ is obviously onto,
the reflection $\refl{u}$ extends to an involutive automorphism of~$\htype$.
\qed

\subsection{Proof of \autoref{prop.A7+A6+A5}}\label{proof.A7+A6+A5}
One has
$\dT=\<\frac78\>\oplus\<\frac67\>\oplus\<\frac43\>
 \oplus\<\frac32\>\oplus\<\frac32\>$, and the image of $\OG(\sset)$ in
$\Aut\dT$ is generated by the reflections $\refl{{\Ga_i}}$, $i=1,2,3$.
Furthermore, one has $\tSigma_2(\bT)=\GG_0^{--}$ and the map
$\det_+\:\Aut\dT\to\{\pm1\}$ is well defined.
Applying \autoref{lem.det+}
with $p=2$,
one finds that
$\det_+\refl{{\Ga_1}}=1$ and $\det_+\refl{{\Ga_2}}=\det_+\refl{{\Ga_3}}=-1$.
In particular, it follows that the homological type is symmetric,
\ie, $\all_\ns(\sset)$ consists of a single real component.

Up to sign, any involutive isometry $a\in\OG(\bT)$ with $\det_+a=-1$ is a
reflection, $a=\pm\refl{x}$ for some $x\in\bT$, $x^2>0$:
one can take for~$x$ a primitive vector generating the
$(-1)$-eigenlattice of $\pm a$, whichever has rank one.
As explained above, $\refl{x}$ must induce $-{\id}$ in one \emph{and only one}
of the components $\dT_3$, $\dT_7$. Hence, $x^2=2^kq$, where $k=1,3$ and
$q=3,7$. (Recall that $x\in(\frac12x^2)\bT\dual$; if $k=2$, then
$\xi:=\frac12x\in\dT_2$ has square $0\bmod\Z$ and $\refl\xi$ is not in the
image of $\OG(\sset)$.)
Obviously, $\eta:=\frac1qx$ is a generator of $\dT_q$; on the other hand, one
can see that $\eta^2/\Ga^2\notin(\Z_q\units)^2$, where $\Ga=\Ga_2$ or~$\Ga_3$
for $q=7$ or~$3$, respectively.
This is a contradiction.
\qed

\subsection{End of the proof of \autoref{th.torus}}\label{s.real.torus}
As in \autoref{s.real.ns}, one can easily see that each set of
singularities~$\sset$ can be obtained by a symmetric perturbation from a
maximizing real one, see \autoref{tab.torus}. Furthermore, the perturbation
can be chosen of \emph{torus type}, \ie, each inner singular point
of weight~$w$
is perturbed to a collection of points of total weight~$w$.
Such perturbations are known to preserve the torus structure.
Hence, by \autoref{cor.real.pert},
the space $\all_3(\sset)$ contains a real curve.
\qed

\subsection{Adjacencies of the strata}\label{s.adjacencies}
Recall that,
with the exception of the set of singularities $\sset=\singset{2A9}$,
the spaces $\all_1(\sset)/\conj$ are connected for all
non-maximizing sextics
(see \autoref{th.classification}).
Together with \cite[Proposition 5.1.1]{degt:8a2} and
\cite{Looijenga:perturbations},
this fact gives us a clear picture of the adjacencies of the
\emph{real} strata; the
only doubtful case of the two components of $\all(\singset{2A9})$
is treated in \autoref{prop.2A9}.

Consider the adjacency graph~$\sC$ of the strata
$\all_\ns(\sset)\subset\all_\ns$ containing
non-real components, and let $\tC$ be the adjacency graph of these non-real
components. One can interpret the vertices and edges of~$\sC$ as,
respectively,
asymmetric primitive homological types and isomorphism classes of
their degenerations,
whereas those of~$\tC$ are oriented homological types and their orientation
preserving degenerations.
With two exceptions, \viz. \singset{A14+A4+A1} and \singset{A13+A6}, see
\autoref{tab.double}, a vertex of~$\sC$ is determined by the corresponding
set of singularities. Most degenerations are of corank one,
in which case a degeneration $\sset'\dg\sset$ is uniquely determined by the
pair $(\sset',\sset)$, see, \eg,~\cite{Dynkin}.
The forgetful projection $\tC\to\sC$ is a double covering, and we are
interested in the structure of this map, in particular, in the connected
components of~$\tC$.

The graph~$\sC$ has several isolated vertices,
\viz. \singset{D7+A10+A2}, \singset{D5+A14},
three vertices representing \singset{A14+A4+A1}, and all
maximizing sets of singularities that are also represented by real curves.
The rest splits into
three larger components, which we denote by $\sC_p$, $p=2,3,7$, and call
\emph{clusters}.
For a fixed~$p$,
the vertices of $\sC_p$ are all sets of singularities in
\autoref{tab.nonreal} containing a $[\,\cdot\,]_p$ pattern and all their
\emph{asymmetric} degenerations, see Figures~\ref{fig.C2}--\ref{fig.C7}.
Denote by $\tC_p\subset\tC$ the pull-back of~$\sC_p$, $p=2,3,7$.
Each double covering
$\tC_p\to\sC_p$ is described by its characteristic class, which we denote by
$\Go_p\in H^1(\sC_p;\FF2)$.

%

Let $\cluster_p:=\bigcup\all_\ns(\sset)$, the union running over all
$\sset\in\sC_p$, $p=2,3,7$. These subspaces of~$\all$ are also called
\emph{clusters}; their connected components are in a one-to-one
correspondence with those of~$\tC_p$.

\figure
\centerline{\cpic{C2}}
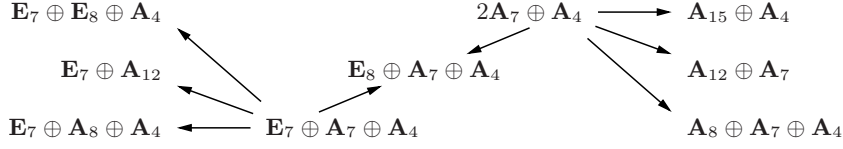
\caption{The graph $\sC_2$}\label{fig.C2}
\endfigure

The graph~$\sC_2$ is shown in \autoref{fig.C2}. Since it is simply connected,
we have the following immediate statement.

\proposition\label{prop.p=2}
The double covering $\tC_2\to\sC_2$ is trivial. Hence,
the cluster
$\cluster_2$ consists of two complex conjugate components.
\done
\endproposition

\figure
\centerline{\cpic{C3}}
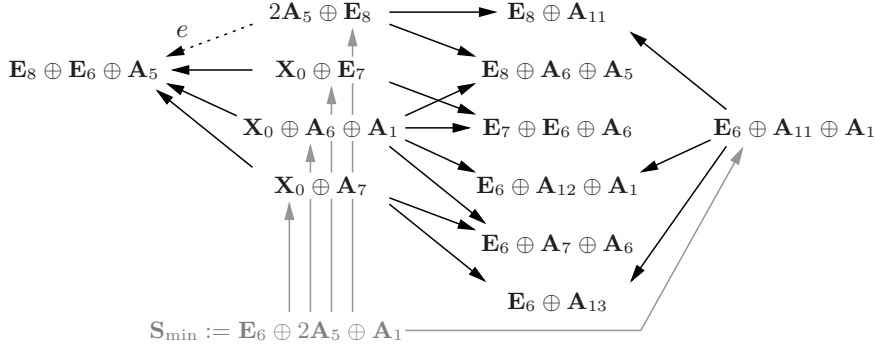
\caption{The graph $\sC_3$ (where $\bX_0:=\bE_6\splus\bA_5$)}\label{fig.C3}
\endfigure

The graph~$\sC_3$ is depicted in \autoref{fig.C3}, where {\em only corank one
degenerations are shown}.
This graph has a minimal vertex $\smin:=\singset{E6+2A5+A1}$, shown in grey.
The closure of $\cluster_3$ contains four real strata
\[*
\singset{2E6+A5+A1}\dg\singset{E7+2E6},\ \singset{2E6+A7},\ \singset{2E6+A6+A1}.
\]
In all four, the real structure interchanges the two $\bE_6$
points; for the non-maximizing set of singularities \singset{2E6+A5+A1}, this
fact can be proved similar to \autoref{lem.asymmetric}.

Regarded as a diagram, $\sC_3$ is not quite commutative.
There are \emph{two} isomorphism classes of degenerations
$\smin\dg\singset{E8+E6+A5}$; in the self-explanatory notation, they are
\[
\singset{[E6+A1]+[A5]+[A5]},\ \singset{[A5+A1]+[E6]+[A5]}\dg\singset{E8+E6+A5}.
\label{eq.smin}
\]
The former factors through the edge
$e\:\singset{2A5+E8}\dg\singset{E8+E6+A5}$ represented by a dotted arrow in
\autoref{fig.C3}, and the latter factors through the three other edges ending
at \singset{E8+E6+A5}. Denote by $e\dual\in H^1(\sC_3;\FF2)$ the class
sending a cycle~$\Ga$, regarded as a sequence of undirected edges, to the
multiplicity of~$e$ in~$\Ga$. Formally, $e\dual$ is the image of the
generator of the group $H^1(e,\partial e;\FF2)=\FF2$
under the relativization homomorphism
$H^1(e,\partial e;\FF2)=H^1(\sC_3,\sC_3\sminus e;\FF2)\to H^1(\sC_3;\FF2)$.

\proposition\label{prop.p=3}
The characteristic class~$\Go_3$ of the double covering $\tC_3\to\sC_3$ is
$\Go_3=e\dual\ne0$. In particular, the cluster $\cluster_3$ is connected.
\endproposition

\proof
Let~$\sC_3'$ be the graph obtained from~$\sC_3$ by removing the (open)
edge~$e$, and let $\tC_3'\subset\tC_3$ be the pull-back of~$\sC_3'$. As
explained above, $\sC_3'$ is a commutative diagram. Hence, the restricted
covering $\tC'_3\to\sC'_3$ is trivial: an orientation of the homological type
extending~$\smin$ induces an orientation of all other homological types.
On the other hand, both degenerations~\eqref{eq.smin} factor through
\singset{2E6+A5+A1} and differ by a transposition of the two $\bE_6$ type points,
which extends to an orientation reversing automorphism of the homological
type. Hence, the double covering $\tC_3\to\sC_3$ is not trivial and the
obstruction is~$e\dual$.
\endproof

The graph~$\sC_7$ is depicted in \autoref{fig.C7}, where shown in black
are the vertices and edges constituting undirected cycles.
(There are two vertices corresponding to the set of singularities
\singset{A13+A6}, see \autoref{tab.double}, each connected by an edge to
\singset{3A6}.)
The group $H_1(\sC_7;\FF2)\cong\FF2^3$ is generated by the three
four-edge cycles~$\Gg_1$, $\Gg_2$, $\Gg_3$,
and the characteristic class $\Go_7$ is
determined by its values on these cycles.

\figure
\centerline{\cpic{C7}}
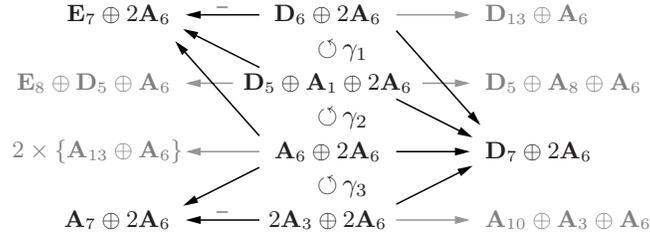
\caption{The graph $\sC_7$}\label{fig.C7}
\endfigure

\proposition\label{prop.p=7}
The characteristic
class~$\Go_7$
of the double covering
$\tC_7\to\sC_7$
is $\Gg_1,\Gg_3\mapsto1$, $\Gg_2\mapsto0$.
In particular, the cluster $\cluster_7$ is connected.
\endproposition

\proof
Consider a quadratic $\FF7$-module~$\CX$.
Recall that the group $\Aut\CX$ is generated by reflections and
there are well defined
homomorphisms $\det,\spin\:\Aut\CX\to\{\pm1\}$
sending a reflection $\refl\xi$ to~$(-1)$ and the class
$14\xi^2\bmod(\Z_7\units)^2\in\Z_7\units/(\Z_7\units)^2=\{\pm1\}$,
respectively,
see,
\eg,~\cite{Cassels}.
Assuming that $\ls|\CX|\cdot\det_7\CX=1\bmod(\Z_7\units)^2$,
define a \emph{$\spin$-orientation} of~$\CX$ as a class of
orthogonal bases $\Ga:=\{\Ga_1,\ldots,\Ga_\ell\}$, $\Ga_i^2=\frac27\bmod2\Z$,
two bases~$\Ga'$, $\Ga''$ being equivalent if the isometry
$\Gs\:\Ga_i'\mapsto\Ga_i''$, $i=1,\ldots,\ell$, has
$\spin\Gs=1$.
Note that the order or the signs of the basis vectors are not important:
isometries reversing the $\spin$-orientation are more subtle.
In particular, the group $\discr_7\sset$ for any $\sset\in\sC_7$ has a
\emph{canonical} $\spin$-orientation.

Let $\ds:=\det\cdot\spin$. In a similar way, using bases with
$\Ga_i^2=-\frac27\bmod2\Z$, we can define the notion of
\emph{$\ds$-orientation} for a $\FF7$-module~$\CY$
satisfying $\ls|\CY|\cdot\det_7\CY=(-1)^\ell\bmod(\Z_7\units)^2$, where
$\ell:=\ell(\CY)$.
An anti-isometry $\CX\to\CY$ takes $\spin$-orientations to
$\ds$-orientations.
There
is a \emph{unique}
$\ds$-orientation on $\<-\frac27\>$;
hence, a $\ds$-orientation on~$\CY$ induces a
$\ds$-orientation on any \emph{codimension one} submodule
$\CZ\subset\CY$ satisfying
$\ls|\CZ|\cdot\det_7\CZ=(-1)^{\ell-1}\bmod(\Z_7\units)^2$.
A similar statement holds for $\spin$-orientations.

The essence of the proof of \autoref{lem.asymmetric} is the fact that,
for any vertex $\sset\in\sC_7$, one has
$\Im[\map_7\:\!\OG_+(\bT)\to\Aut\dT_7]\subset\Ker\ds$.
(If $\Gm(\sset)=19$, this follows from~\cite{Shimada:maximal}.)
Hence, there is a bijection $\conv\:\pss\mapsto\dso$ between positive sign
structures on~$\bT$ and $\ds$-orientations on~$\dT_7$.
(The particular choice of $\conv$ is not important; it can be
fixed separately for each isomorphism class.)
Thus, an oriented homological type $(\htype,\pss)$
can be declared \emph{positive} or
\emph{negative} according to whether
the anti-isometry $\dS\to\dT$
does or does not take the
canonical $\spin$-orientation of~$\dS_7$ to $\conv(\pss)$.

Given a lattice extension $\iota\:\sset\fe\sset'$, the homomorphisms
$\iota\otimes\Q$ and $\iota\dual$ induce additive relations
$\iota_*\:\dS_7\dashrightarrow\dS'_7$ and
$\iota\dual\:\dS'_7\dashrightarrow\dS_7$.
If $\iota$ is one of the \emph{black} arrows in \autoref{fig.C7}, both
$\iota_*$ and $\iota\dual$ are true homomorphisms; they give rise,
in a canonical way, to either
an isomorphism $\dS_7=\dS_7'$ or a splitting
$\dS_7=\dS_7'\oplus\<\frac27\>$ (if $\sset=3\bA_6$),
which respect the canonical $\spin$-orientation. Passing to the
transcendental lattices, we conclude that,
in either case, a $\ds$-orientation on~$\dT_7$
induces one on~$\dT_7'$. On the other hand, $\bT'\subset\bT$ is a maximal
positive definite sublattice and $\bT$ and $\bT'$ have a \emph{common}
positive sign structure $\pss=\pss'$. Hence, we can assign to~$\iota$ a sign
$\epsilon=\pm1$ so that the $\ds$-orientation on~$\dT'_7$ induced by
$\conv(\pss)$ equals $\epsilon\conv(\pss')$.
This sign depends on the conventions~$\conv$, but the product
$\epsilon:=\epsilon_1\epsilon_2\epsilon_3\epsilon_4$ over a four-edge cycle
$c:=(\iota_1,\iota_2,\iota_3,\iota_4)$ does not, as each convention is used
twice.
It is immediate from the
definitions that $\epsilon=(-1)^{\Go_7(c)}$.
Now, the statement of the proposition is proved by a routine computation of
the signs, \cf. \autoref{ex.C7} below.
\endproof

\example\label{ex.C7}
We illustrate the computation of the signs in the previous proof. All
rank two lattices involved are of the form $\Z u\oplus\Z v$,
$u^2=2^{r+1}\cdot7$, $v^2=2^{s+1}\cdot7$, $r,s\ge0$, and \emph{for such
lattices}, we define $\conv$ to take the positive basis $\{u,v\}$ to a
basis $\{\Ga_1,\Ga_2\}$ with $\Ga_1:=\frac17(2^{r+2}u+2^{s+1}v)$.
(For a module of length two, one vector of square $-\frac27\bmod2\Z$ is
enough to define a $\ds$-orientation.
For the comparison purposes, it is convenient to consider the basis
$\Gb_1:=\frac17\cdot2^ru$, $\Gb_2:=\frac17\cdot2^sv$ with
$\Gb_1^2=\Gb_2^2=\frac27\bmod2\Z$, so that $\Ga_1=4\Gb_1+2\Gb_2$.
In terms of the $\Gb$-basis, the transposition of the two vectors or changing
the sign of one of them reverses the $\ds$-orientation.)
To avoid choices for rank three lattices, we consider a pair of arrows
$\sset'\gd\sset\dg\sset''$.
Let $\sset=\singset{D6+2A6}$ (the
topmost
pair in \autoref{fig.C7}). Then $\bT=\Z u\oplus\Z v\oplus\Z w$, $u^2=v^2=14$,
$w^2=-2$, and $\bT',\bT''\subset\bT$ are spanned, respectively, by
$u':=u$, $v':=v$ and $u'':=3u+7w$, $v'':=v$.
(We use the fact that each transcendental lattice involved is
known to be unique in its
genus and merely `guess' a representation producing the correct discriminant.
Since we know that the sign is well defined, it suffices to consider a
particular pair of sublattices.)
The orientations of the two bases
are coherent,
and the coefficient $3\notin(\Z_7\units)^2$ in the expression
for~$u''$ tells us that the product of the signs associated with this pair of
arrows is $(-1)$: one has $\Gb''_1=-\Gb'_1$ and $\Gb''_2=\Gb'_2$.
To complete
the cycle $\gamma_1$, consider the other rank three lattice
$\bar\sset=\singset{D5+A1+2A6}$ and its degenerations
$\sset'\gd\bar\sset\dg\sset''$
(the second
pair in
\autoref{fig.C7}).
Then, in the self-explanatory notation,
one has $\bar\bT=\Z\bar u\oplus\Z\bar v\oplus\Z\bar w$,
$\bar u^2=14$,
$\bar v^2=28$, $\bar w^2=-2$,
and the generators of $\bT',\bT''\subset\bar\bT$ can be
represented as $u':=\bar u$, $v':=2\bar v+7\bar w$ and
$u'':=\bar u$, $v'':=\bar v$.
Since $2\in(\Z_7\units)^2$, we can deduce that
$\bar\Gb_1=\Gb'_1$ and $\bar\Gb_2=\Gb'_2$.
Hence, the cumulative sign of the cycle~$\Gg_1$ is
$\epsilon=1\cdot (-1)\cdot 1\cdot 1=-1$,
\ie, $\Go_7(\Gg_1)=1$.

A similar computation, slightly more involved if $\sset=3\bA_6$
(for which one has $\bT=(\Z u+\Z v)\oplus\Z w$,
$u^2=v^2=0$, $u\cdot v=7$, $w^2=14$),
shows that
the sign convention for rank three lattices can be chosen so that only
the two
arrows marked with a `$-$' in \autoref{fig.C7} have associated sign~$(-1)$;
this proves \autoref{prop.p=7}.
\endexample

\let\.\DOTaccent
\def\cprime{$'$}
\bibliographystyle{amsplain}
\bibliography{degt}

\end{document}